\documentclass{article}
\usepackage{mathtools}
\usepackage{mathrsfs,enumerate}
\usepackage{amsmath}
\usepackage{amssymb}
\usepackage[numbers]{natbib}
\usepackage{hyperref}
\usepackage{dsfont}
\usepackage[draft]{todonotes}
\usepackage[a4paper,top=2cm,bottom=2cm,left=2cm,right=2cm,bindingoffset=5mm]{geometry}

\numberwithin{equation}{section}

\usepackage{color}
\usepackage{dsfont}
\usepackage[latin1]{inputenc}  
\textwidth16cm \numberwithin{equation}{section} \vspace{8cm}
\usepackage{hyperref}

\newtheorem{theorem}{Theorem}[section]
\newtheorem{proposition}[theorem]{Proposition}
\newtheorem{corollary}[theorem]{Corollary}
\newtheorem{lemma}[theorem]{Lemma}
\newtheorem{definition}[theorem]{Definition}
\newtheorem{remark}[theorem]{Remark}

\newcommand{\R}{{\mathbb R}}

\newcommand{\cL}{{\mathcal L}}
\newcommand{\cM}{{\mathcal M}}
\newcommand{\cP}{{\mathcal P}}

\def\cW{{\mathcal W}}
\newcommand{\vep}{\varepsilon}

\newcommand{\vfi}{\varphi}

\def\lg{\langle}
\def\rg{\rangle}
\def\tv{\tilde v}
\def\ty{\tilde y}

\def\L1k{L^1(\lg x\rg^{k})}

\def\mdep{m_{\de,\vep}}

\newcommand{\be}{\begin{equation}}
\newcommand{\ee}{\end{equation}}
\newcommand{\rife}[1]{{(\ref{#1})}}

\newcommand{\dive}{{\rm div}}

\newcommand{\into}{{\int_{\Omega}}}

\newcommand{\intd}{{\int_{\R^{2d}}}}

\def\de{\delta}

\def\vep{\varepsilon}

\def\vfi{\varphi}
\def\qed{{\unskip\nobreak\hfil\penalty50
          \hskip2em\hbox{}\nobreak\hfil\mbox{\rule{1ex}{1ex} \qquad}
   \parfillskip=0pt
   \finalhyphendemerits=0\par\medskip}}

\def\rife#1{(\ref{#1})}

\title{Long-time contractivity estimates for kinetic Kolmogorov-Fokker-Planck equations}

\author{
Nicol\`o Forcillo\thanks{Michigan State University, Department of Mathematics, 619 Red Cedar Road, East Lansing, MI 48824.  Email: \texttt{forcill1@msu.edu}. "Partially supported by  the MIUR Excellence Department Project  Math@TOV" awarded to the Department of Mathematics, University of Rome Tor Vergata. The author was also partially supported by the ``Problemi variazionali/nonvariazionali: interazione tra metodi integrali e principi del massimo" 2023-2024 INDAM-GNAMPA project.} 
		\and Alessio Porretta\thanks{Department of Mathematics, University of Rome Tor Vergata. 
		Via della Ricerca Scientifica 1, 00133 Roma, Italy.    Email: \texttt{porretta@mat.uniroma2.it}. "Partially supported by  the Excellence Project MatMod@TOV of the Department of Mathematics of the University of Rome Tor Vergata and Italian (EU Next Gen) PRIN project 2022W58BJ5 ({\it PDEs and optimal control methods in mean field games, population dynamics and multi-agent models}, CUP E53D23005910006). The author is also a member of a GNAMPA research group of Indam. 
	}  
 }

\begin{document}

\maketitle
\begin{abstract} 
We prove long-time contractivity estimates and exponential rates of convergence to equilibrium for solutions of hypoelliptic diffusion equations, which include the well-known Kolmogorov equation and similar kinetic Fokker-Planck equations in $\R^d$. Compared to the existing literature, our proof exploits a different approach, elementary and self-contained, based on oscillation estimates for the adjoint problem. We first prove contractivity in Wasserstein distances through doubling variables (coupling) methods. Next, we upgrade the estimate to weighted $L^1$-(or total variation) norms,  thanks to short-time hypocoercivity gradient estimates.  
\end{abstract}

%
\vskip1em

\tableofcontents


\section{Introduction}

This paper concerns the long-time convergence of solutions  of kinetic-type Fokker-Planck equations, such as 
\be\label{KFP}
\begin{cases}
	\partial_t m - \Delta_v m + \dive_y ( H(v,y) m) =  \dive_v ( B(v,y)  m )\,,  \qquad (v,y)\in \R^{2d}\,, &  \\
	m(0,v,y)=m_0(v,y)\,, & 
\end{cases}
\ee
where $m_0$ is a probability measure on $\R^{2d}$, and $B,H:\R^{2d}\to \R^d$ are functions such that $D_vH$ is positive definite. The simplest and most famous example is given by the Kolmogorov equation (corresponding to $H(v,y)=v$) 
\be\label{KO}
	\partial_t m - \Delta_v m + v \cdot \nabla_y m =  \dive_v ( B(v,y)  m )\,,
\ee
which is the prototype of degenerate convection-diffusion  equations that exhibit a hypoelliptic behavior. We recall that the solution of \rife{KO} describes the law of the Langevin process 
$$
	\begin{cases}
		dy_t = v_t \, dt,  & \\
		dv_t= - B(v_t, y_t) \, dt + \sqrt 2\, dW_t\,, & 
	\end{cases}
$$
where $W_t$ is a Brownian motion in $\R^d$.
It is well known that solutions of \rife{KO} converge exponentially fast,  as $t$ goes to infinity,  towards the stationary state,  if either $y$ is in a compact setting (e.g., on the torus) or $B$ is a suitable confining vector field.  
The occurrence of this exponential stability is a subtle combination of diffusion and transport, and it is strictly related to the underlying hypoellipticity of the operator. 

There exist many different methods to show that exponential stabilization occurs. Originally, in the  model case, the explicit knowledge of the stationary state and/or of the fundamental solution was used to derive this kind of result. However, an approach that avoids the use of explicit solutions appears preferable for tackling more general problems with a possibly similar structure. In \cite{Villani}, Villani introduced the notion of hypocoercivity to frame the question of long-time stability, showing that  an exponential contractivity estimate, such as
$$
\|{m_1(t)}-{m_2(t)}\|_X\leq C e^{-\omega t} \|m_{01}-m_{02}\|_X\,,
$$
could be proved in the Hilbert setting of suitable energy spaces, being $X$, typically, a weighted $H^1-$ or $L^2-$space. This is also related to spectral gap methods, already used before to quantify  regularizing effects and long-time rate of convergence, see, e.g., \cite{Desvillettes+V}, \cite{H-Nier}, \cite{herau}, \cite{herau-nier}. Later on, those approaches have been refined to a large extent and merged with entropy methods, geometric Poincar\'e inequalities, and different tools from semigroup theory (\cite{Arnold-Erb}, \cite{Bakry+al}, \cite{Dolbeault+al}, \cite{Duan}, \cite{Gualdani+al}, \cite{MiMou}, \cite{zhang}).    
Different techniques, more directly based on Harris' recurrence theorem for Markov semigroups (\cite{HM}), or the so-called Meyn-Tweedy approach (\cite{MeynTweed}, \cite{tweed}), were employed in other works, see, for instance, \cite{Canizo+al}, \cite{Cao}, \cite{mattingly+al}. Their  connection with our present  article will be further commented below.

The goal of this paper is to introduce a new method to prove the exponential rate of convergence,  extending to degenerate Kolmogorov operators the techniques recently developed in  \cite{Po} for convection-diffusion equations with uniform or fractional diffusion. This approach derives the exponential decay by duality from oscillation estimates of the adjoint  problem. Precisely, the focus is moved from the conservative equation \rife{KFP} to the dual advection-diffusion problem
\be\label{AD}
\partial_t u - \Delta_v u  -   H(v,y)  \cdot \nabla_y u +  B(v,y) \cdot \nabla_v u=0\,,
\ee
where the  decay  of  $u$ is proved in suitably weighted seminorms { (see estimate \rife{ktheta})}.  
%
%
%

This idea has its roots in the coupling methods developed in probability (\cite{ChenLi},  \cite{Lin-Rog}) and, in parallel, in doubling variables techniques of the PDE community (\cite{CIL}, \cite{IL}), in the same spirit as in \cite{PoPr}. The outcome is a self-contained, elementary  approach, where pointwise oscillation estimates for the advection-diffusion equation \rife{AD} lead to the exponential decay of \rife{KFP}  bypassing any spectral analysis, Poincar\'e inequality, or integral Harnack inequality. Somehow, our work answers a natural question (mentioned by Villani in \cite[Chapter 2]{villani_book}), namely, whether the exponential stabilization of the (hypoelliptic) kinetic Fokker-Planck equation could be proved through a simple coupling method, as in the diffusive case. A positive answer was already given in  \cite{Eberle+al} by probabilistic arguments, showing  the exponential contractivity  in  Wasserstein distances for solutions of \rife{KO}. Here, we complement and extend  \cite{Eberle+al} from a  PDE perspective. Our contribution is twofold. First, we provide an entirely PDE proof of the decay in Wasserstein distances, only based on  maximum principle, which gives similar conclusions for the general problem \rife{KFP} as those obtained in \cite{Eberle+al} for the Kolmogorov equation. Next, we improve the result to achieve the decay in total variation norm, by combining the oscillation estimates with short-time $L^\infty \to W^{1,\infty}-$smoothing effect. This latter tool seems interesting in its own, and should be considered as the $L^\infty-$version, on the adjoint problem, of the $L^2 \to H^1-$smoothing effect lying at the heart of Villani's hypocoercivity estimate.  
The final outcome can be summarized in the following contractivity estimate, when the initial data are probability measures with finite first moments, denoted by $\cP_1(\R^{2d})$.

\begin{theorem}\label{main}  Assume that $H,B$ satisfy  \rife{Blip} and \rife{H1}-\rife{H2}, and there exists a Lyapunov function $\vfi$ (cf. Definition \ref{def-Lyap}) fulfilling conditions \rife{liplocvfi}, \rife{vfi2}, and \rife{vfi3}. 
Then, there exist $\omega, K>0$ such that, for every initial data $m_{01}, m_{02}\in \cP_1(\R^{2d})$ for which  $\vfi\in L^1(dm_{0i})$, $i=1,2$, the corresponding solutions $m_1,m_2$ of \rife{KFP} satisfy 
\be\label{output}  
\| {m_1(t)}-{m_2(t)}\|_{TV_\vfi} \leq K \, e^{-\omega t} \, \|m_{01} - m_{02}\|_{TV_\vfi}, 
\ee
where $\|\cdot\|_{TV_\vfi}$ is the norm of the total variation weighted by $\vfi$.
\end{theorem}

Let us spend a few words on the hypotheses of the above result.  

A crucial role in Theorem \ref{main} is played by the assumption that $D_vH$ is positive definite (cf. \rife{H1}). Roughly speaking, this means that the $v$-linearization of the drift field $H$ is nondegenerate. The case in which $D_vH$ is negative definite would, of course, be the same.  
When $H$ is linear with constant coefficients, our assumptions can be reduced to one form of the H\"ormander condition (though not the most general one), in particular, the H\"ormander condition of first order (say, only using first-order commutators). We discuss this issue in Remark \ref{hypo-rem}. It has to be mentioned that, in the special case of linear drifts with constant coefficients, the exponential convergence under the general H\"ormander condition is proved  in \cite{Arnold-Erb} by using entropy methods and the explicit form of the steady state. We also hope to extend our approach to incorporate the case of general H\"ormander operators. So far, our assumptions  include  and generalize  the case of the Kolmogorov equation, appearing as the  first step in the direction of treating general nonlinear and inhomogeneous  first-order terms. 

{ The requirement, in Theorem \ref{main}, that a suitable Lyapunov function exists should not be surprising, since the state space is noncompact and  thus some confinement condition is necessary for long-time stabilization. This is a classical ingredient used to show the ergodicity property of the Markov process associated with the operator \rife{AD}. Explicit examples of  Lyapunov functions are given at the end of the paper, precisely in Section \ref{examples}.  Similar forms of these functions already appear in  previous works (see, e.g., \cite{Bakry+al}, \cite{Canizo+al}, \cite{Cao}, \cite{Eberle+al}, \cite{mattingly+al}) with  either analytic or probabilistic approaches. Some of the aforementioned papers  rely on the version of Harris' theorem proved in  \cite{HM}, whose cornerstone is indeed the existence of a Lyapunov function. We stress again, however, that our method departs from this literature because we bypass the use of Harnack inequalities for the Fokker-Planck equation in favor of global oscillation estimates proved on the adjoint problem \rife{AD}.  In this context, the contractivity estimate is reduced to a maximum principle argument on the adjoint problem, where coupling methods naturally come into play}.

Let us remark that the estimate of Theorem \ref{main} holds for   weak (distributional) solutions of \rife{KFP}. The key point indeed is to show that weak solutions are dual to viscosity solutions of the adjoint problem \rife{AD}, which are sufficiently robust for the main estimate. This duality, already exploited in \cite{Po} for other drift-diffusion operators, is also very natural to handle sign-changing solutions. In this direction, the estimate \rife{output} can alternatively be reformulated as the decay of zero-average solutions, see Theorem \ref{zero-average}.  
%
%
 
As we mentioned before, in the intermediate steps needed to prove Theorem \ref{main}, we also establish other results that may have interest  in their own. They can be found in the next sections, so it is worth giving a quick plan of the paper. In Section \ref{setting}, we present some notation and  the setting of our assumptions, and we discuss the notion of solution to \rife{KFP} and \rife{AD} we are going to use. We also state an existence and uniqueness result (Theorem \ref{teo-ex}) for both the adjoint problems,  which is the natural background of our estimates. Section \ref{oscillation} is devoted to prove  the main estimate for problem \rife{AD}, namely, the exponential decay of suitable oscillation seminorms, see Theorem \ref{LTO}. An immediate byproduct is the decay in Wasserstein distances of solutions of \rife{KFP}, see Theorem \ref{decay-Wass}.
In Section \ref{hypo}, we show the $L^\infty\to W^{1,\infty}-$smoothing effect for the equation \rife{AD}, by suitably rephrasing Villani's hypocoercivity argument to obtain pointwise gradient bounds through the classical Bernstein method, see  Proposition \ref{hypocoer}. A weighted version, including the Lyapunov function growth into the estimate, is also needed to our purposes.
Finally, in Section \ref{final-decay}, we prove our main results, Theorem  \ref{main-new} (which includes the above Theorem \ref{main}) and Theorem \ref{zero-average}, and we give some examples of application, where confinement properties of the drift ensure the existence of Lyapunov functions with suitable conditions. 

\section{Setting of the problem}\label{setting}

\subsection{Notation}

Let us start with some notation. In the following, $\R^d$, $d\geq 1$, denotes the $d$-dimensional Euclidean space, where $dy$ represents the Lebesgue measure. $I_d$ is the {$d$-dimensional} identity matrix, and $|y|$ is the usual Euclidean norm, for $y\in \R^d$. Given two real numbers $a,b$, we use the standard notation $ a\vee b= \max(a,b)$ and $a\wedge b= \min(a,b)$. 
For vectors $z,w\in \R^d$, $z\otimes w$ is the matrix such that $(z\otimes w)_{ij}= z_i w_j$.  For any {$d\times d$ (real)} matrix $A$, tr$(A)$ denotes its trace.

The space $C_b(U)$ is the set of continuous bounded functions on an open  set $U$, and $C_c(U)$ is the set of continuous compactly supported ones. 
The space $\cM(\R^d)$ denotes the set of finite (signed) Borel measures $m$, whose  total variation norm is $\| m\|_{TV}{:=} |m|(\R^d)$. Here, $|m|= m^+ + m^-$, where $m^\pm$ are positive measures giving the Hahn decomposition of $m$. { The subset of probability measures is represented by $\cP(\R^d)$, endowed with the topology of narrow convergence. This is  equivalent to the weak-$*$ convergence of measures in $\cM(\R^d)$, induced by duality with the uniform convergence in   $C_b(\R^d)$,  see, e.g., \cite[Chapter 5]{AGS}.  }We refer to $\cM_k(\R^d)$ as the subset of $\cM(\R^d)$ consisting of measures $m$ with finite $k$th moment, i.e., $|v|^k\in L^1(\R^d, d|m|)$. In particular, the set $\cP_1(\R^d)$ denotes the space of probability measures with finite first moment, where we consider the Kantorovich-Rubinstein distance
\begin{equation}\label{Was-dist}
\cW_1(m,\tilde m):= \sup\bigg\{\int_{\R^{d}} \zeta \, d(m-\tilde m):\,\, \zeta\,\hbox{ Lipschitz,}\,\,  \|D\zeta\|_\infty\leq 1 \bigg\}\,.
\end{equation}

We fix a positive, continuous function $\vfi$ on $\R^d$. $L^1(\vfi(x))$ represents the usual Lebesgue space $L^1(\R^d, \mu)$, defined in terms of the measure $\mu = \vfi(x) dx$, and $L^\infty(\vfi(x))$ is the space of measurable functions $u:\R^d\to \R$ such that $u\, \vfi $ is essentially bounded (namely, $u\vfi \in L^\infty(\R^d)$).  
Similarly, we denote by $\cM_\vfi(\R^d)$ the subset of measures $m$ such that $\vfi\in L^1(\R^d, d|m|)$, and we set  ${ \|m\|_{TV_\vfi}}{:=} \int_{\R^{d}} \vfi\, d|m|$.

For given $T>0$, we set $Q_T:= (0,T)\times \R^d$. We denote by $C^{1,2}(Q_T)$ the space of functions, defined on $Q_T$, that are $C^1$ in $t$ and $C^2$ in $x$.  We refer to $C([0,T];\cM(\R^d)^*)$ as the space of functions $m$  from $[0,T]$ into $\cM(\R^d)$, which are continuous with respect to the weak$^*$-convergence of measures.  { Similarly, we use $C([0,T];\cP(\R^d))$ when working with probability measures, and $C([0,T];\cP_1(\R^d) )$ for continuous functions $m$ from $[0,T]$ into $\cP_1(\R^d)$, endowed with the $\cW_1$-distance}.

\subsection{Assumptions { and well-posedness of} the differential problems}

Let $v,y\in \R^d$ and denote by $\nabla_v, {\dive_v}, \Delta_v$, etc...,  the corresponding differential operators in the $v$-variable (the same for $y$). Given $\kappa \geq 0$, we consider the problem  
\be\label{FP}
	\begin{cases}
		\partial_t m - \Delta_v m - \kappa \Delta_y m + \dive_y ( H(v,y) m) =  \dive_v ( B(v,y)  m )\,,   & t\in (0,T)\,, \quad  (v,y)\in \R^{d}\times\R^d\,,\\
		m(0,v,y)=m_0(v,y)\,, & (v,y)\in \R^{d}\times\R^d\,,
	\end{cases}
\ee
where $m_0 \in \cM(\R^{2d})$. We recall that this equation is associated to  the stochastic system 
$$
	\begin{cases}
		dy_t = H(v_t, y_t)\,dt  + \kappa \sqrt 2\,dW_t & \\
		dv_t = - B(v_t, y_t)\,dt + \sqrt 2\, d\tilde W_t 
	\end{cases}
$$
for two independent $d$-dimensional Brownian motions $W_t,   \tilde W_t$, since ${m(t)}$ can be  interpreted as the law of the  process $(v_t, y_t)$ in $\R^{2d}$, see, e.g., \cite{Anceschi+al}. 
The case in which $y$ is taken in the flat torus could also be considered, as well as the case $v\in \R^d, y\in \R^k$ with $k\neq d$. We will comment those two possible extensions later, {see Remark \ref{extensions}}.  Here, we focus on the symmetric case, which includes the Langevin equation of Brownian particles when $\kappa=0$, $H=v$, and, formally,  $\ddot{y_t} = -B(\dot{y_t}, y_t)  + \sqrt 2  \frac{dB_t}{dt} $.

We introduce the linear operator 
\be\label{def-lin-oper}
	\cL[\vfi]: = -\Delta_v \vfi  - \kappa \Delta_y \vfi -  H(v,y)  \cdot \nabla_y \vfi +  B(v,y) \cdot \nabla_v \vfi\,,
\ee
which defines the adjoint problem to \rife{FP}. Indeed, we will also consider this adjoint PDE (put in forward {time} direction), which is the following 
\be\label{pb}
	\begin{cases}
		\partial_t u - \Delta_v u  - \kappa\,\Delta_y u -  H(v,y)  \cdot \nabla_y u +  B(v,y) \cdot \nabla_v u=0\,,   & t\in (0,T)\,,\quad  (v,y)\in \R^{d}\times\R^d\,,
		\\
		u(0,v,y)=u_0(v,y)\,, & (v,y)\in \R^{d}\times\R^d\,.
	\end{cases}
\ee
In \rife{FP} and \rife{pb}, we assume that $B:\R^d\times \R^d \to \R^d$ satisfies the following  Lipschitz condition  
\be\label{Blip}
|B(v_1,y_1)-B(v_2, y_2)| \leq \ell_B (|v_1-v_2|+ |y_1-y_2|)\,, 
\ee
for some $\ell_B>0$. About  the function $H: \R^d\times \R^d \to \R^d$, we assume that it holds
\begin{eqnarray}
	& D_vH(v,y)\geq \gamma\, I_d \qquad \qquad \qquad  \forall (v,y)\in \R^{2d}\,, & \label{H1} \\
	 & |H(v_1,y_1)-H(v_2, y_2)| \leq \ell_H (|v_1-v_2|+ |y_1-y_2|) \qquad \forall (v_1,y_1), (v_2,y_2)\in \R^{2d}\,,  & \label{H2}
\end{eqnarray}
with $\gamma>0$, $\ell_H>0$. Note, incidentally, that \rife{H1} and \rife{H2} imply $\ell_H\geq \gamma$.

\begin{remark}\label{hypo-rem}
The nondegeneracy condition \rife{H1} plays a crucial role in what follows, ensuring the hypoellipticity of the drift-diffusion operator.

Let us discuss this condition in the linear case, i.e., if $H(v,y)= Cv+ Dy$,  where $C,D$ are constant matrices in $\R^d$.
There,  assumption \rife{H1} corresponds to requiring $C$ to be a {positive} definite matrix. To fix the ideas, we assumed $C>0$, but, of course, the case $C<0$ is equally treated up to changing $y$ into $-y$ in the equation. We also stress that if $\kappa=0$, then the change of variable  $My=  y' $ would turn equation \rife{FP} into a similar {one} where $C$ is replaced by $MC$. So, actually, {we only need} the {condition} $MC>0$ for some invertible matrix $M$. This is the case, for instance, if $C$ is just assumed to be invertible, rather than positive definite.

To sum up, when $\kappa=0$ and the drift is linear, our assumptions are equivalent to requiring the so-called first-order H\"ormander condition for linear drift-diffusion operators $-{\rm Tr}(AD^2(\cdot))+ FD(\cdot)$, namely, that no subspace of  $Ker(A)$ is mapped into $Ker(A)$ by $F^*$ (in our setting, if $H(v,y)= Cv+ Dy$, then this would mean  $Ker(C)\neq 0$).
As pointed out, for instance, in \cite[Section 2]{Arnold-Erb},  the full H\"ormander condition is actually more general, only asking that no subspace of $Ker(A)$ is left-invariant by $F^*$. However, this more general condition requires to use commutators of order higher than one. We expect to exploit this further generality in the future. By the way, we stress  
once more that our setting is neither confined to linear drifts nor to constant coefficients, obviously.  
\end{remark}

Let us now specify the notion of solution we use for the above initial value problems. 
As for \rife{FP}, we work with the standard notion of weak (distributional) solution.

\begin{definition}\label{weak-sol} 
Let $m_0\in \cM(\R^{2d})$. A function $m\in {C([0,T];\cM(\R^{2d})^*)}$  is a weak solution to \rife{FP} if
$$
 \int_0^T\intd  \left( -\partial_t \phi+ \cL[\phi] \right) \,d{m(\tau)} \,d\tau = \intd {\phi(0,v,y)}\, dm_0 
 \qquad\forall \,  \phi\in C^{1,2}_c([0,T) \times \R^{2d} )\,.
 $$
 \end{definition}

\vskip0.3em

Concerning problem \rife{pb}, which is a degenerate second-order problem, it seems very natural to exploit the general concept of viscosity solutions, nowadays classical because of \cite{CIL}. This notion preserves the maximum principle (the only ingredient needed for the estimates) and   is stable for vanishing viscosity (i.e., as $\kappa\to0$ above). 
As a byproduct of the estimates we prove later, we get the following existence and uniqueness result, { whose proof is postponed to the Appendix.}

\begin{theorem}\label{teo-ex} Assume that \rife{Blip} and \rife{H1}-\rife{H2} hold, and $\kappa\geq 0$. For any $T>0$, we have:
\begin{itemize}

\item[(i)] Given $u_0 \in C_b(\R^{2d})$, there exists a unique viscosity solution { $u\in C_b(\overline Q_T)$ of  \rife{pb}. Moreover,    for every $t>0$, {$u(t)$} is Lipschitz continuous and satisfies estimate \rife{hypo}.}

\item[(ii)] Given { $m_0\in \cP_1(\R^{2d})$,  there exists a unique weak solution $m\in C([0,T];\cP(\R^{2d}) )$  of  \rife{FP}.  In addition, $m\in C([0,T];\cP_1(\R^{2d}) )$ and
\be\label{dual-form}
\begin{split}
	&
	 \intd \zeta\, d{m(t)}  = \intd {u(0,v,y)}\, dm_0 
	 \\
	 & \qquad \qquad\forall \, t\in (0,T)\,,\, \zeta\in C_b(\R^{2d}),  \text{ and }u\in C_b([0,t] \times \R^{2d})\quad \mbox{such that} \\
	 &  \qquad\qquad\qquad \hbox{$u$   is  a viscosity solution of }\,\,  \mbox{\bf (B)}\,\,\begin{cases}
	-\partial_t u+ \cL[u]     =  0 &  \text{ in } ( 0,t ) \times \R^{2d}\,, \\ 
	u(t,v,y) = \zeta (v,y)\qquad  &  \text{ in }   \R^{2d}\,.
\end{cases} 
\end{split}
\ee
For general (possibly signed) measures $m_0 \in \cM_1(R^{2d})$,  there exists a unique weak solution $m$ of  \rife{FP} such that $ \sup_{t\in (0,T)} \intd |z|\, d{|m|(t)}  <\infty$, and in addition $m$ satisfies \rife{dual-form}.
}
\end{itemize}

\end{theorem}

{ 
\begin{remark} There is a slight difference in the above statement concerning the uniqueness of solutions of \rife{FP}, depending on whether ${m(t)}$ is a probability measure or a general (possibly signed) measure. In this latter case, we need to require a priori that ${|m|(t)}$ has finite first moments, bounded in $(0,T)$. By contrast, for positive measures, this estimate can be deduced directly from the formulation. This subtle distinction between probability solutions (nonnegative and mass-preserving) and general solutions (possibly signed and not a priori mass-preserving) is classical for Fokker-Planck equations, see, e.g.\cite{boga+al}. However, the two notions coincide here, provided that the first moments are controlled.
\end{remark}

\vskip0.3em

Let us comment further on the statements of Theorem \ref{teo-ex}. The well-posedness of viscosity solutions of \rife{pb} appears as a variant of the general theory. However, we could not find a result directly applied to \rife{pb}, so we will give a proof here for completeness. Uniqueness of merely {\it weak solutions} to \rife{FP} is not trivial in general, although it is not surprising here due to the Lipschitz growth of the drift terms. For nondegenerate diffusions, this is well known, see, for instance, \cite{boga+al}, \cite{shapo}. For the Kolmogorov equation, a similar result can be found in \cite{Anceschi+al}. Consequently, the most interesting outcome of the previous theorem is  the  duality equality \rife{dual-form} between weak and viscosity solutions of the two adjoint problems, as well as the consistency of the problems with the vanishing viscosity limit. As observed in \cite{Po}, the dual representation \rife{dual-form} can eventually serve itself as a definition of robust solutions of Fokker-Planck equations, whenever the ordinary theory does not apply anymore. As a matter of fact, the formulation \rife{dual-form} suits possibly signed measures, and it guarantees positivity and mass preservation. A statement in this sense is the following corollary. 

\begin{corollary}\label{cor-ex}  Let $m_0\in \cM_k({\R^{2d}})$ for some $k>0$. Then, there exists a unique  $m\in C([0,T];\cM({\R^{2d}})^*)$ satisfying \rife{dual-form}. In addition, $m_0\geq 0$ implies ${m(t)}\geq 0$ for all $t\in (0,T]$, and $\intd d{m(t)}= \intd dm_0$ for every $t>0$.
\end{corollary}

}

\section{Oscillation decay estimates}\label{oscillation}

Together with the hypoellipticity condition \rife{H1}, a crucial role will also be played by Lyapunov functions since we are going to investigate the long-time stabilization for a noncompact state space. 

\begin{definition}\label{def-Lyap} We say that $\vfi(v,y)$ is a Lyapunov function for  the operator $\cL$ if  it satisfies
\be\label{lyapu}
	\vfi\in C^2(\R^d\times \R^d)\,, \,\, \vfi(v,y)\mathop{\longrightarrow}^{|(v,y)|\to \infty} \infty\,,\hbox{ and }\,\, \exists \,\, \omega_0>0\,:\quad \liminf_{|(v,y)|\to \infty} \, \frac{\cL[\vfi]}{\vfi}\geq \omega_0\,.
\ee
Without loss of generality, we will also suppose $\vfi\geq 1$ in $\R^d\times \R^d$.
\end{definition}

Let us mention that a condition such as \rife{lyapu} is satisfied under suitable coercivity assumptions on the drift terms, for which we refer to Subsection \ref{examples}.

In connection with the Lyapunov function $\vfi$, we will consider { (viscosity) solutions $u$ of \rife{pb} so that $u =o(\vfi)$ as $|(v,y)|\to \infty$ (uniformly in $t\in (0,T)$), namely,
\be\label{uvsvfi}
\lim_{|(v,y)|\to \infty}\, \left[\sup_{t\in (0,T)}\, \frac{u(t,v,y)}{\vfi(v,y)} \right] =0 \,.
\ee
Of course, this condition is trivially fulfilled for bounded solutions, but it seems convenient (and costless!) to consider the broader setting of possibly unbounded solutions, like, for instance, unbounded Lipschitz solutions.}

\subsection{Exponential  decay of  oscillation seminorms}

In this subsection, we prove the exponential decay for a weighted  oscillation seminorm of the solution. For $\theta\in (0,1]$ and $w\in C(\R^{2d})\cap L^\infty(\vfi^{-1})$, we define
\be\label{def-seminorm}
	[w]_\theta: = \sup_{(v,y)\neq  (\tv, \ty)} \,\, \frac{|w(v,y)- w(\tv,\ty)|}{(\vfi(v,y)+ \vfi(\tv, \ty))( d^\theta \wedge 1)}\,,\qquad d=| (v,y)- (\tv, \ty)|\,.
\ee
Then, we have the following result.

\begin{theorem}\label{LTO} { Let $u$ be a viscosity solution of \rife{pb} such that $u(t)$ is Lipschitz continuous in $\R^{2d}$ (with a uniform Lipschitz constant  for  $t\in (0,T)$) and \rife{uvsvfi} is true, for some  Lyapunov function $\vfi$ associated to the operator $\cL$}. 
Assume that \rife{Blip} and \rife{H1}-\rife{H2} hold, and, for some $\ell_\vfi>0$,
\be\label{liplocvfi}
	| \nabla_v \vfi( v,y)-\nabla_v\vfi(\tv,  \ty)|\leq \ell_\vfi\, [\vfi( v,y)+ \vfi( \tv,\ty)] (|v-\tv| + |y-\ty|) \,.
\ee
Then, for every $\theta\in (0,1]$,  there exist $\omega, K>0$ such that 
\be\label{ktheta}
	[{u(t)}]_\theta\leq K\, e^{-\omega t} [u_0]_\theta 
\ee
for every $t>0$, where $\omega, K$ depend  on $\gamma, \ell_H, \ell_B, \vfi,$ and $\theta$.
\end{theorem}

\begin{proof} 
We first give the proof for any $\theta<1$. In the end, we explain the modifications needed to handle the limiting case $\theta=1$.

We introduce the rotated norm
\be\label{rotated-norm}
	|| (v,y) || := |v+\mu y |+ \lambda |y|\,, \qquad \hbox{where $\mu=  \frac{2\ell_H}{\gamma}$\,, $\quad \lambda=\frac4{\gamma} ((\mu+1) \ell_H+ 2\ell_B)$\,.}
\ee
Note that  $\mu$ and $\lambda$ only depend on the constants $\gamma, \ell_H,$ and $\ell_B$ appearing in \rife{Blip}, \rife{H1}, and \rife{H2}. { The precise value of $\mu$ and $\lambda$ is not especially relevant and could have been fixed later in the proof. However, we observe, since now, that  $\mu\geq 2$ and $\lambda >4(\mu+1)$ being $\ell_H\geq \gamma$, as noticed before. 

Of course, the norm \rife{rotated-norm} is equivalent to the standard Euclidean metric in $\R^{2d}$. In particular, since $\lambda \geq 1+\mu $,   we have }
$$
|v-\tilde{v}|+|y-\tilde{y}| \leq |v-\tilde{v}+\mu(y-\tilde{y})|+(1+ \mu) |y-\tilde{y}|\le  |v-\tilde{v}+\mu(y-\tilde{y})|+\lambda {|y-\tilde{y}|}{  = || (v-\tv,y-\ty)||}\,.
$$
We regularize  $|| \cdot ||$ as 
	$$
	|| (v,y) ||_\vep:= |v+\mu y |_\vep+ \lambda |y|_\vep\,,\quad \hbox{with  $|z|_\vep: = \sqrt{\vep^2+ |z|^2}-\vep$\,.}
	$$
For later convenience, we denote by
\be\label{def-xi}
\xi:= v+\mu y\,.
\ee
Next,  we define the function $\psi$ as
\be\label{def-psi}
	\psi(r):= 1-e^{-C_2 r^\theta}\,,
\ee
for $\theta\in (0,1)$ and some $C_2\geq 1$ to be determined later. We also fix a generic time horizon $T>0$.  We claim that, for some $C_2, \omega, K,$ and $L$ to be  chosen later, we have
\be\label{claim}
 	u(t,v,y)- u(t, \tv, \ty ) \leq  e^{-\omega t}  \left(K[\vfi(v,y)+ \vfi(\tv,\ty)]+ L\right) \psi(\rho) + \frac{\eta}{T-t}\,,
\ee 
for every $t\in (0,T), v,y,\tv,\ty\in \R^d$ and $\eta$ sufficiently small, where 
$$
	\rho{=} || (v-\tv,y-\ty)||\,.
$$ 
In order to prove \rife{claim}, we argue by contradiction. We define 
$$
	{\mathbf m}:= \sup_{t, (v,y), (\tv,\ty)} \bigg\{ u(t,v,y)- u(t, \tv, \ty ) -   e^{-\omega t}  \left(K[\vfi(v,y)+ \vfi(\tv,\ty)]+ L\right) \psi(\rho) -  \frac{\eta}{T-t} \bigg\}   
$$
and we suppose that ${\mathbf m}>0$. We are going to show that this leads to a contradiction up to a suitable choice of  $K, L, C_2$ (sufficiently large) and $\omega$ sufficiently small.

We first notice that ${\mathbf  m}$ is finite and is a global maximum attained at some point. Precisely, we have 
\be\label{psibass}
	\psi(\rho)\geq \frac12 (C_2\rho^\theta \wedge 1) \geq \frac12 (C_2 \rho  \wedge 1)
\ee
since $1-e^{-t}\geq \frac12 (t\wedge 1)$. 
Being, by assumption, $u(t,v,y)= o(\vfi(v,y))$ as $|(v,y)|\to \infty$ (uniformly in $(0,T)$),  it holds $|u(t,v,y)|\leq  e^{-\omega t}  \frac K4 \vfi( v,y)+C_{K,\omega}$. In addition, $u(t,v,y)- u(t, \tv, \ty )\leq \ell^T_u\, \rho$ for some $\ell_u^T>0$,  due to the Lipschitz continuity and {the equivalence between $|| \cdot ||$ and the standard Euclidean metric in $\R^{2d}$}. Hence, we can estimate
$$ 
	u(t,v,y)- u(t, \tv, \ty )    \leq   e^{-\omega t}  \left( \frac K4 [\vfi( v,y)+ \vfi( \tv,\ty)] + \tilde C_T\right) (C_2 \rho  \wedge 1)\,,
$$
for some constant $\tilde C_T$ (possibly depending on $u$), which yields, by \rife{psibass},
\begin{align*}
	u(t,v,y)- u(t, \tv, \ty )  & -   e^{-\omega t}  \left(K[\vfi(v,y)+ \vfi(\tv,\ty)]+ L\right) \psi(\rho)  \\ 
	& \leq   e^{-\omega t}  \left( -\frac  K2 [\vfi( v,y)+ \vfi( \tv,\ty)] + 2\tilde C_T-L\right) \psi(\rho)\,.
\end{align*}
Since $\vfi$ blows-up at infinity, we deduce  that
\be\label{limsup-inf}
	\limsup_{|(v,y)| \vee |(\tv, \ty)| \to \infty} \,\, \bigg\{ u(t,v,y)- u(t, \tv, \ty ) -   e^{-\omega t}  \left(K[\vfi( v,y)+ \vfi( \tv,\ty)]+ L\right) \psi(\rho) -  \frac{\eta}{T-t} \bigg\} \leq 0\,,
\ee
uniformly for $t\in (0,T)$.  Therefore, assuming   ${\mathbf m}>0$ implies (by continuity of $u$)  that the sup is attained at some point in $[0,T]\times \R^{2d}$.

Now, for $\vep>0$ small, we consider the regularized distance $\rho_\vep=|| (v-\tv,y-\ty)||_\vep$ and, correspondingly, the value 
\be\label{defmeps}
	{\mathbf m}_\vep:= \sup_{t, (v,y), (\tv,\ty)} \bigg\{ u(t,v,y)- u(t, \tv, \ty ) -   e^{-\omega t}  \left(K[{\vfi(v,y)+ \vfi(\tv,\ty)}]+ L\right) \psi(\rho_\vep) -  \frac{\eta}{T-t} \bigg\}   \,.
\ee
Notice  that $|z|_\vep\leq |z|$, and so
\be\label{rhotoeps}
	\rho_\vep\leq \rho\leq  \rho_\vep + (1+\lambda)  \vep\,.
\ee
In particular,  being $\psi$ increasing, we have ${\mathbf m}_\vep\geq {\mathbf m}>0$. 
Moreover, we claim  that \rife{limsup-inf} even holds with $\rho$ replaced by $\rho_\vep$. Specifically, repeating the argument exploited above, for any small $\sigma>0$, we obtain that there exists a constant $C_\sigma$ such that  
$$
	u(t,v,y)- u(t, \tv, \ty )    \leq 
	[\sigma (\vfi(v,y)+\vfi(\tv,\ty)) + C_\sigma] \mathds{1}_{\{C_2\rho_\vep^\theta> 1\}} + \ell^T_u \rho\, \mathds{1}_{\{C_2\rho_\vep^\theta\leq 1\}}  \,,
$$
for some $\ell_u^T>0$. 
So, from \rife{rhotoeps}, we infer (recalling that $C_2\geq 1$ and $\theta\in (0,1)$)
\begin{align*}
	u(t,v,y)- u(t, \tv, \ty )     & \leq [\sigma (\vfi(v,y)+\vfi(\tv,\ty)) + C_\sigma] \mathds{1}_{\{C_2 \rho_\vep^\theta> 1\}} +  \ell^T_u\,     \rho_\vep \mathds{1}_{\{C_2\rho_\vep^\theta\leq 1\}} + \ell^T_u (1+\lambda)  \vep \\
	& \leq  [\sigma (\vfi(v,y)+\vfi(\tv,\ty)) + C_\sigma] \mathds{1}_{\{C_2 \rho_\vep^\theta> 1\}} +  \ell^T_u\,  C_2    \rho_\vep^\theta  \mathds{1}_{\{C_2\rho_\vep^\theta\leq 1\}} + \ell^T_u (1+\lambda)  \vep
	\\ 
	& \leq [ \sigma (\vfi(v,y)+\vfi(\tv,\ty)) + C_\sigma + \ell^T_u ] 2\, \psi(\rho_\vep)+ \ell^T_u (1+\lambda)  \vep
	\,,
\end{align*}
where we have also used \rife{psibass} in the last step.
Hence, we get, for any $t\in (0,T)$,  
\begin{align*}
	&u(t,v,y)- u(t, \tv, \ty ) -   e^{-\omega t}  \left(K[\vfi(v,y)+ \vfi(\tv,\ty)]+ L\right) \psi(\rho_\vep) -  \frac{\eta}{T-t}\\  
	&\leq  \left\{ [C_\sigma + \ell_u^T ] 2-   e^{-\omega T}  \left( ( K- 2\sigma e^{\omega T}) [\vfi(v,y)+ \vfi(\tv,\ty)]+ L\right)\right\} \psi(\rho_\vep) + \ell^T_u (1+\lambda)  \vep-  \frac{\eta}{T-t}\,.
\end{align*}
Taking $\sigma<e^{-\omega T}K/2$ and $\vep$ such that $ \vep \ell^T_u (1+\lambda)  < \frac\eta T$, we achieve \rife{limsup-inf} with $\rho$ replaced by $\rho_\vep$. Consequently, ${\mathbf m}_\vep$ is a global maximum as well, attained at some point  $(t_\vep,v_\vep,y_\vep,\tv_\vep,\ty_\vep)$. 
We also notice that, obviously,  $t_\vep <T$ and $t_\vep>0$, provided   
\be\label{u0}
	u_0(v,y)- u_0( \tv, \ty ) \leq      \left(K[\vfi( v,y)+ \vfi( \tv,\ty)]+ L\right) \psi(\rho_\vep)+ \frac\eta T
\ee
for every $v,y,\tv,\ty$ and $\vep $ sufficiently small.  This is the case, for instance,  if it holds (recall \rife{def-seminorm})
\be\label{Ku0}
	K \geq  2[u_0]_\theta\,.
\ee
Indeed, due to \rife{rhotoeps}, if $\rho \geq 2(1+\lambda) \vep$, then we have $\rho\leq 2\, \rho_\vep$. Thus, using the Lipschitz character of  $u_0$ and \rife{psibass}, and taking $C_2$ large, it holds
\begin{align*}
	u_0(v,y)- u_0( \tv, \ty )  & \leq  [u_0]_\theta\, [\vfi( v,y)+ \vfi( \tv,\ty)] ((2\rho_\vep)^\theta \wedge 1) + \ell_{u_0} \, \rho\, \mathds{1}_{\{ \rho<2(1+\lambda) \vep\}} 
	\\ & \leq  [u_0]_\theta\, [\vfi( v,y)+ \vfi( \tv,\ty)] ((C_2\rho_\vep)^\theta \wedge 1) + \ell_{u_0} \, 2(1+\lambda)\vep
	\\ & \leq  [u_0]_\theta\, [\vfi( v,y)+ \vfi( \tv,\ty)] 2 \psi(\rho_\vep) + \ell_{u_0} \, 2(1+\lambda)\vep\,.
\end{align*}
Therefore, \rife{u0} is satisfied as long as $K\geq 2[u_0]_\theta$ and $2(1+\lambda)\vep\, \ell_{u_0}\leq \frac\eta T$. We conclude that $(t_\vep,v_\vep,y_\vep,\tv_\vep,\ty_\vep)$ is a local (and global) maximum point, where ${\mathbf m}_\vep$ is attained.
\vskip1em

By classical viscosity solutions theory (\cite{CIL}), this implies that there exist matrices $V,  \tilde V, Y, \tilde Y$ such that  
\be\label{matri}
	\left(
	    \begin{array}{cccc}
	      V  & \ldots & 0  & 0 \\
	      \ldots & Y & 0 & 0  \\
	       0 & 0 & - \tilde V & \ldots \\
	      0 &0 & \ldots & -\tilde Y \\
	    \end{array}
  	\right)
 \leq    e^{-\omega t} D^2_{(v,y), (\tv,\ty)}\zeta\,,
\ee
where
$$
	\zeta:=   \left(K[\vfi( v,y)+ \vfi( \tv,\ty)]+ L\right) \psi(\rho_\vep)\,, 
$$
and
\be\label{zeta1}
\begin{split}
	& \frac\eta{(T-t)^2} - \omega \,e^{-\omega t} \zeta  -  e^{-\omega t } \left(H(v,y) \cdot D_y \zeta +  H(\tilde v, \ty) \cdot  D_{\tilde y}\zeta\right) 
	\\
	& \quad + e^{-\omega t}   \left(  B(v,y)\cdot \nabla_v \zeta + B(\tv,\ty)  \cdot \nabla_{\tilde v} \zeta\right)
	\leq   {\rm tr}(  V -  \tilde V )   + \kappa \,  {\rm tr}(  Y -  \tilde Y )\,. 
\end{split}
\ee
Here, we have dropped the index $\vep$ from the points $(v,y),(\tv,\ty)$ to lighten the notation, and so we will do henceforth.

To handle the diffusion terms,  we introduce the coupling matrix   
{ \be\label{couple}
	\left(
		\begin{array}{cccc}
		  I_d  &0 & {\mathcal C} & 0 \\
		 0 & \kappa I_d & 0 & \kappa I_d  \\
		  {\mathcal C}^* & 0 & I_d&0 \\
		  0 & \kappa I_d & 0 & \kappa I_d \\
		\end{array}
  	\right)\,,
\ee}
which we assume to be nonnegative. We also suppose that ${\mathcal C}$ is symmetric to ease computations. { Note that we are only using the synchronous coupling in the $y$-variable, since in the spirit of our result, the $y$-diffusion is meant to be vanishing}.

Multiplying \rife{matri} by the above matrix and taking the traces, we get
$$
	{\rm tr}(   V -   \tilde V )+ \kappa \,  {\rm tr}(  Y -  \tilde Y ) \leq e^{-\omega t} {\rm tr}\left(  D^2_{vv}\zeta+ 2 {\mathcal C} D_{\tv x}^2\zeta+   D^2_{\tv\tv}\zeta\right)+ \kappa \,e^{-\omega t} {\rm tr}\left(  D^2_{yy}\zeta+ 2  D_{\ty y}^2\zeta+   D^2_{\ty\ty}\zeta\right)\,, 
$$ 
so that \rife{zeta1} implies
\begin{align*}
	& \frac\eta{(T-t)^2} - \omega \,e^{-\omega t} \zeta  -  e^{-\omega t } \left(H(v,y) \cdot D_y \zeta +  H(\tilde v, \ty) \cdot  D_{\tilde y}\zeta\right) 
	  + e^{-\omega t}   \left(  B(v,y) \cdot \nabla_v \zeta + B(\tv,\ty)  \cdot \nabla_{\tilde v} \zeta\right)
	\\
	& \qquad \leq  e^{-\omega t} {\rm tr}\left(   D^2_{vv}\zeta+ 2 {\mathcal C}D_{\tv v}^2\zeta+  D^2_{\tv\tv}\zeta\right) + \kappa \,e^{-\omega t} {\rm tr}\left(  D^2_{yy}\zeta+ 2  D_{\ty y}^2\zeta+   D^2_{\ty\ty}\zeta\right) \,.
\end{align*}
Using the explicit form of the test function $\zeta$, dropping the first term, and dividing by $e^{-\omega t}$, it follows that, in view of \rife{def-lin-oper},
\be\label{zeta2}
	\begin{split}
		&  \psi(\rho_\vep) K  \left\{  \cL[\vfi](v,y)- \omega \vfi(v,y) + \cL[\vfi](\tv, \ty) - \omega\vfi(\tv,\ty) \right\} - \omega \,L \psi(\rho_\vep)   { + }  \left(K[\vfi( v,y)+ \vfi( \tv,\ty)]+ L\right) \psi'(\rho_\vep) 
		\\ 
		&  \times \left( -(H(v,y)-H(\tv,\ty)) \cdot \nabla_y \rho_\vep  + ( B(v,y) -B(\tv,\ty) ) \cdot \nabla_v \rho_\vep \right)  \leq   \left(K[\vfi( v,y)+ \vfi( \tv,\ty)]+ L\right)
		\\
		&\times 2 {\rm tr}\left( I_d-   {\mathcal C} \right)  D^2_{vv}(\psi(\rho_\vep))+  
		2 K \psi'(\rho_\vep)  \left( I_d-   {\mathcal C} \right)\nabla_v\rho_\vep \cdot \left( \nabla\vfi( v,y)-\nabla  \vfi(\tv, \ty)\right)  \,.
	\end{split}
\ee
Let us now specify our choice of ${\mathcal C}$. For $\xi$ given by  \rife{def-xi} and $\tau\geq 0$, we take
\be\label{choice-coup}
	{\mathcal C} = I- \tau \frac{\xi-\tilde \xi}{|\xi-\tilde \xi|} \otimes \frac{\xi-\tilde \xi}{|\xi-\tilde \xi|}\,.
\ee
With this, we observe that the matrix in  \rife{couple} is nonnegative provided $0\leq \tau\leq 2$. This selection of coupling includes what is called the synchronous coupling ($\tau=0$) and the reflection one ($\tau=2$), and we will choose one of these two options depending on the distance $\rho_\vep$ between the points $(v,y)$ and $(\tv,\ty)$.  Notice, indeed, that we need to take $\tau=0$ whenever $\xi= \tilde \xi$. Otherwise, ${\mathcal C}$ would be singular.

Let us further elaborate the second-order terms. According to the definition of $\rho_\vep$, we compute 
\begin{align*}
	& D^2_{vv} \psi(\rho_\vep)   = \psi''(\rho_\vep) \nabla_v\rho_\vep\otimes \nabla_v\rho_\vep + \psi'(\rho_\vep) D^2_{vv}\rho_\vep
	\\ & = \psi''(\rho_\vep)  \bigg(\frac{\xi-\tilde \xi}{|\xi-\tilde \xi|_\vep+\vep} \otimes \frac{\xi-\tilde \xi}{|\xi-\tilde \xi|_\vep + \vep} \bigg) + \psi'(\rho_\vep) \frac1{|\xi-\tilde \xi|_\vep + \vep} \bigg(I_d- \frac{\xi-\tilde \xi}{|\xi-\tilde \xi|_\vep + \vep} \otimes \frac{\xi-\tilde \xi}{|\xi-\tilde \xi|_\vep + \vep} \bigg)\,.
\end{align*}
Therefore, we obtain
\be\label{second-der}
	2 {\rm tr} \left(\left( I_d- {\mathcal C} \right)  D^2_{vv}(\psi(\rho_\vep))\right) =
	2\tau \, \psi''(\rho_\vep)\, \left( \frac{|\xi-\tilde \xi|}{|\xi-\tilde \xi|_\vep + \vep}\right)^2 
	+ 2  \tau \,\frac{\psi'(\rho_\vep)}{|\xi-\tilde \xi|_\vep + \vep}\, \frac{\vep^2}{(|\xi-\tilde \xi|_\vep + \vep)^2}\,.
\ee
We now show that  \rife{zeta2}, together with \rife{second-der}, leads to a contradiction, with suitable choices of the parameters.  We analyze three different cases, depending on the distance between the points $(v,y)$ and $(\tv,\ty)$.
\vskip0.4em
{\it Case 1 ($\rho$ large)}: We first analyze the case that $\rho_\vep \geq R_1$, for some large $R_1>1$ to be fixed later.

We estimate, {by virtue of \rife{Blip} and \rife{H2}}, since  $|\nabla_v \rho_\vep| \leq 1, |\nabla_y \rho_\vep| \leq (\mu+\lambda)$, and recalling that $|v-\tilde{v}|+|y-\tilde{y}|\le \rho$,
\be\label{HB-est}
	\begin{split}
		& | (H(v,y)-H(\tv, \ty)) \cdot \nabla_y \rho_\vep| +    | ( B(v,y) -B(\tv,\ty) )  \cdot \nabla_v \rho_\vep| 
		 \\ & \qquad \leq  ( \ell_H (\mu+ \lambda) + \ell_B) ( |v-\tv|+ |y-\ty|) \leq  ( \ell_H (\mu+ \lambda) + \ell_B) (\rho_\vep + (1+\lambda) \vep) \leq  C \rho_\vep\,,
	\end{split}
\ee
where we have also exploited \rife{rhotoeps}  and the fact that $\rho_\vep$ is large enough. Here, $C$ depends on  $\mu, \lambda, \ell_B$, and $\ell_H$ (hence, it only depends on $\ell_B,\ell_H,$ and $\gamma$).
Next, we take $\tau=0$ in the choice of the coupling {\rife{choice-coup}} (i.e., $\mathcal C= I_d$). Hence, we obtain, from \rife{zeta2},
\be\label{aprile}
	\begin{split}
		&  \psi(\rho_\vep) K  \left\{  \cL[\vfi](v,y)- \omega \vfi(v,y) + \cL[\vfi](\tv, \ty) - \omega\vfi(\tv,\ty) \right\} - \omega \,L \psi(\rho_\vep) 
		\\ 
		&   -   \left(K[\vfi( v,y)+ \vfi( \tv,\ty)]+ L\right)  C  \psi'(\rho_\vep)\rho_\vep   \leq 0\,.
	\end{split}
\ee 
Since  \rife{def-psi} implies 
\be\label{psi'-times-r}
	\psi'(r)r=\theta C_2  r^\theta e^{-C_2r^\theta}\,,
\ee
it results that $\psi'(r)r$ is decreasing for  $C_2 r^\theta\geq 1$, so, a fortiori, if $r\ge 1$. Then, we have $\psi'(\rho_\vep)\rho_\vep\leq \psi'(R_1)R_1$ for  $\rho_\vep \ge R_1\geq 1$. Being $C_2\ge 1$, 
it also holds  $1/2 \leq \psi(\rho_\vep)\leq 1$  for $\rho_\vep \ge R_1$.  Consequently, we get
$$
 C \psi'(\rho_\vep)\rho_\vep  \leq 2C \psi(\rho_\vep) \psi'(R_1)R_1\,.
$$
Renaming  $\tilde C= 2C$, it follows, in view of \rife{aprile},
\be\label{preKL}
	\begin{split}
		&  \psi(\rho_\vep)\Big\{ K \big(   \cL [\vfi](v,y)- (\omega+ \tilde C \psi'(R_1)R_1)   \vfi(v,y) + \cL [\vfi](\tv, \ty) -  (\omega+ \tilde C\psi'(R_1)R_1) \vfi(\tv,\ty)\big)   \\ 
		&  \qquad   - L( \omega + \tilde C\psi'(R_1)R_1)\Big\}
		\leq   0\,.
	\end{split}
\ee 
We recall that {$|z|_\vep\leq |z|$ and}
$$
	\rho_\vep \leq || (v,y)|| + || (\tv,\ty)||\,.
$$
%
Therefore,  $\rho_\vep\geq R_1$ implies that  $ || (v,y)||\geq  R_1/2$ or $|| (\tv,\ty)|| \geq  R_1/ 2 $. Without loss of generality, we assume the first holds. Up to choosing $\omega$ small and $R_1$ large so that
\be\label{omegar1}
	\omega+ \tilde C\psi'(R_1)R_1<\frac{\omega_0}2\,,
\ee
where $\omega_0$ is  given by \rife{lyapu}, we have  that $\cL[\vfi]- (\omega+ \tilde C\psi'(R_1)R_1) \vfi\geq -k_0$ for some $k_0>0$. Hence, we estimate, according to \rife{lyapu} and  \rife{omegar1} again,
\begin{align*}
	\cL[\vfi](v,y)- (\omega+ \tilde C\psi'(R_1)R_1)   \vfi(v,y)  & + \cL[\vfi](\tv, \ty) -  (\omega+ \tilde C\psi'(R_1)R_1) \vfi(\tv,\ty) 
	\\ & \geq \cL[\vfi](v,y)- \frac{\omega_0}2 \vfi(v,y)-k_0 
	 \geq \frac{\omega_0}4 \vfi(v,y)\,, 
\end{align*}
as long as  $R_1$ is sufficiently large, only depending on $\omega_0$ and $\vfi$ (and $\mu$ and $\lambda$).   
We thus deduce, from \rife{preKL}, that if $\rho_\vep \geq R_1$, then
$$
K \frac{\omega_0}4\, \bigg(\inf_{\{||(v,y)||\geq R_1/2\}} \vfi\bigg)  \leq L( \omega + \tilde  C\psi'(R_1)R_1) \,.
$$
However, this is impossible if we choose $K$ as
\be\label{choiK}
	K= L\, \frac{8( \omega + \tilde C\psi'(R_1)R_1) }{\omega_0 \vfi_{R_1}}\,,\qquad \vfi_{R_1}:= \inf\{\vfi(v,y):\,\, ||(v,y)||\geq R_1/2\}\,.
\ee
 {\it Case 2 ($|\xi-\tilde \xi | > \sqrt \vep$)}  We now discuss the case  $\rho_\vep < R_1$ and $|\xi-\tilde \xi | > \sqrt \vep$. If this occurs, then we have, for $\vep$ small, 
$$
	\rho_\vep \geq |\xi-\tilde \xi |_\vep {\geq |\xi-\tilde \xi |-\vep} \geq \frac12 \sqrt \vep\,, 
$$
thus,  from \rife{rhotoeps},
\be\label{rho_eps-bigger-sum-x-y}
	|v-\tv|+ |y-\ty| \leq \rho \leq \rho_\vep +  (1+\lambda) \vep \leq  C\, \rho_\vep \,.
\ee
Then, we can estimate the drift terms in \rife{zeta2} exactly as in \rife{HB-est}, using again the Lipschitz character of  $H$ and $B$. Similarly, we treat the last term in \rife{zeta2} using  $|\nabla_v\rho_\vep|\leq 1$ and assumption \rife{liplocvfi}, which yields
$$
2 K \psi'(\rho_\vep) ( I_d- {\mathcal C})\nabla_v\rho_\vep \cdot \left(  \nabla \vfi( v,y) -\nabla \vfi(\tv, \ty)\right)  
\leq C \, K    \, [\vfi( v,y)+ \vfi( \tv,\ty)]  \psi'(\rho_\vep)\rho_\vep\,.
$$
Hereafter, we denote by $C$ possibly different constants  only depending on $\gamma, \ell_H, \ell_B$, and $\ell_\vfi$.
Finally, we use \rife{second-der} with $\tau =2 $ (note that $\tau=1$  would equally work) and $ \cL[\vfi]-  \omega  \vfi\geq -k_0$. Overall, we obtain from  \rife{zeta2},
\be\label{nick}
	\begin{split}
		&  -2 \psi(\rho_\vep)k_0\,  K    - \omega \,L \psi(\rho_\vep) -  \left(K[\vfi( v,y)+ \vfi( \tv,\ty)]+ L\right)  C\,  \psi'(\rho_\vep) \rho_\vep  \\
		& \qquad \leq   \left(K[\vfi( v,y)+ \vfi( \tv,\ty)]+ L\right)  
		\bigg\{ 4     \, \psi''(\rho_\vep)\, \bigg( \frac{|\xi-\tilde \xi|}{|\xi-\tilde \xi|_\vep + \vep}\bigg)^2 
		+ 4     \,\frac{\psi'(\rho_\vep)}{|\xi-\tilde \xi|_\vep + \vep}\, \frac{\vep^2}{(|\xi-\tilde \xi|_\vep + \vep)^2}\bigg\}\,.
	\end{split}
\ee 
We now observe that, due to \rife{choiK}, there exists $c>0$ (independent of $C_2, R_1,$ and $L$) such that  
\be\label{stimaKpsi}
	K\psi(\rho_\vep) \leq  c \, L\, \left(  \frac{\omega \psi(\rho_\vep) + \psi'(\rho_\vep)\rho_\vep} { \vfi_{R_1}} \right)\quad {  \forall \rho_\vep<R_1}\,.
\ee
Indeed, if $C_2\rho_\vep^\theta \geq 1$, then { $\psi'(\rho)\rho$ is decreasing for $\rho\geq \rho_\vep$}, so \rife{choiK} implies, for $\rho_\vep <R_1$,  
$$
	K\psi(\rho_\vep) { \leq L\, \frac{8( \omega \psi(\rho_\vep)+ \tilde C\psi'(\rho_\vep)\rho_\vep) }{\omega_0 \vfi_{R_1}}}\,,
$$
which yields \rife{stimaKpsi}. Alternatively,  { we have $C_2\rho_\vep^\theta < 1$}, and then $\psi(\rho_\vep) \leq \frac e\theta \psi'(\rho_\vep)\rho_\vep$ (using \rife{psi'-times-r} and that $1-e^{-t} \leq e\, t e^{-t}$ for $t\leq 1$). Hence, $K\psi(\rho_\vep) \leq c\,  K\, \psi'(\rho_\vep)\rho_\vep$ and \rife{stimaKpsi}  immediately follows because $K\leq  \frac 4{ \vfi_{R_1}} L$, according to \rife{choiK} and \rife{omegar1}. Therefore, \rife{stimaKpsi} holds for every $\rho_\vep <R_1$ and we have
$$
	2 \psi(\rho_\vep)k_0\,  K   + \omega \,L \psi(\rho_\vep) \leq \frac {2k_0c} { \vfi_{R_1}}  \psi'(\rho_\vep)\rho_\vep\, L  + \left( 1+ \frac {2k_0c} { \vfi_{R_1}}\right) \omega \,L \psi(\rho_\vep)\,.  
$$
Note that $\vfi_{R_1}$ is uniformly bounded below (since $\vfi\geq 1$), and eventually it can be taken large up to increasing the value of $R_1$. In particular, we will assume henceforth that $R_1$ is sufficiently large so that 
$$
	\frac {2k_0c} { \vfi_{R_1}} \leq \frac\gamma4\,,
$$
where we recall that $\gamma$ is given in \rife{H1}.
With this in hand, we refine the above estimate as
\be\label{quelk0}
	2 \psi(\rho_\vep)k_0\,  K   + \omega \,L \psi(\rho_\vep) \leq \frac\gamma 4 \psi'(\rho_\vep)\rho_\vep\, L  +  \left( 1+ \frac \gamma 4\right)  \omega \,L \psi(\rho_\vep)\,
\ee
whenever $\rho_\vep< R_1$. 
Inserting this into \rife{nick}, we obtain
\be\label{zeta5}
	\begin{split}
		&  \left(K[\vfi( v,y)+ \vfi( \tv,\ty)]+ L\right)  
		\left\{ - 4      \psi''(\rho_\vep)\, \left( \frac{|\xi-\tilde \xi|}{|\xi-\tilde \xi|_\vep + \vep}\right)^2 
		 -   \frac{4\psi'(\rho_\vep)}{|\xi-\tilde \xi|_\vep + \vep}\, \frac{\vep^2}{(|\xi-\tilde \xi|_\vep + \vep)^2} - C\,\psi'(\rho_\vep) \rho_\vep \right\}  
		\\ & \qquad \leq \left( 1+ \frac \gamma 4\right) \omega \,L \psi(\rho_\vep)\,. 
	\end{split} 
\ee
From the definition of the function $\psi$, it holds
\be\label{zeta6}
	\begin{split}
		& \left\{    - 4     \psi''(\rho_\vep)\, \left( \frac{|\xi-\tilde \xi|}{|\xi-\tilde \xi|_\vep + \vep}\right)^2 
		 -  \frac{4\psi'(\rho_\vep)}{|\xi-\tilde \xi|_\vep + \vep}\, \frac{\vep^2}{(|\xi-\tilde \xi|_\vep + \vep)^2} - C\,\psi'(\rho_\vep) \rho_\vep \right\} 
		 \\ & \quad = C_2 \theta\, \rho_\vep^{\theta-2} e^{-C_2 \rho_\vep^\theta} \left( 4  [ (1-\theta)+ C_2\theta  \rho_\vep^\theta]\left( \frac{|\xi-\tilde \xi|}{|\xi-\tilde \xi|_\vep + \vep}\right)^2 - 4   \rho_\vep \, \frac{\vep^2}{(|\xi-\tilde \xi|_\vep + \vep)^3}  - C \rho_\vep^2\right)
		 \\ & \quad \geq C_2 \theta\, \rho_\vep^{\theta-2} e^{-C_2 \rho_\vep^\theta} \left( 2    (1-\theta)+  2 C_2\theta  \rho_\vep^\theta    -  4  \sqrt \vep  \rho_\vep  - C \rho_\vep^2\right)\,,
	\end{split}
\ee
where we have used  that $|\xi-\tilde \xi |>\sqrt \vep$ gives
\be\label{rooteps}
	\left( \frac{|\xi-\tilde \xi|}{|\xi-\tilde \xi|_\vep + \vep}\right)^2\geq 1-\vep\geq \frac12 \quad \hbox{and} \quad  \frac{\vep^2}{(|\xi-\tilde \xi|_\vep + \vep)^3}\leq \sqrt \vep\,.
\ee
Now, we fix $C_2$, in terms of $R_1$, by imposing   
$$
 2C_2\theta  r^\theta   -  4  \sqrt \vep \, r  - C r^2 \geq 0 \qquad \forall  r\leq R_1\,.
$$
Thus, we deduce, from \rife{zeta5} and \rife{zeta6}, 
$$ 
	\left(K[\vfi( v,y)+ \vfi( \tv,\ty)]+ L\right)   C_2 \theta\, \rho_\vep^{\theta-2} e^{-C_2 \rho_\vep^\theta}  2   (1-\theta)    \leq \left( 1+ \frac \gamma 4\right) \omega \,L \psi(\rho_\vep)\,, 
$$
which yields, since $\rho_\vep\leq R_1$, $\theta\in (0,1)$, and $\psi\leq 1$,
$$
	C_2 \theta\, R_1^{\theta-2} e^{-C_2 R_1^\theta}  2   (1-\theta)    \leq \left( 1+ \frac \gamma 4\right) \omega\,.
$$
Nevertheless, this is impossible if we choose $\omega$ sufficiently small, namely, 
\be\label{omega1}
	\omega < \left( 1+ \frac \gamma 4\right)^{-1} C_2 \theta\, R_1^{\theta-2} e^{-C_2 R_1^\theta}  2   (1-\theta) \,.
\ee
{\it Case 3 ($|\xi-\tilde \xi | \le \sqrt \vep$)} By virtue of the above choices and conditions for the parameters, we are only left with the possibility that $\rho_\vep<R_1$ and  $|\xi-\tilde \xi | \le \sqrt \vep$. This scenario corresponds to $|\xi-\tilde \xi|$ near zero (or possibly zero), where we cannot use the coupling by reflection for the variable $\xi$.  We indeed take ${\mathcal C}=I_d$. Moreover, we exploit $ \cL[\vfi]-  \omega  \vfi\geq -k_0$. Hence, we achieve, according to \rife{zeta2}, 
\be\label{zeta3}
	\begin{split}
		&  -2 \psi(\rho_\vep)k_0\,  K    - \omega \,L \psi(\rho_\vep) 
		 +   \left(K[\vfi( v,y)+ \vfi( \tv,\ty)]+ L\right)  \psi'(\rho_\vep)  \left\{-\left( H(v,y)-H(\tv,\ty)\right) \cdot \nabla_y \rho_\vep \right.  
		\\
		& \left. + ( B(v,y) -B(\tv,\ty) ) \cdot \nabla_v \rho_\vep)\right\} \leq 0 \,.
	\end{split}
\ee
As for the drift terms, recalling that $x= \xi- \mu y$ and $ |\nabla_y \rho_\vep|\leq (\mu+ \lambda)$, we have
\begin{align*}
	& - (H(v,y)-H(\tv,\ty)) \cdot \nabla_y \rho_\vep    = -  \int_0^1 D_v H(sv+ (1-s)\tv, sy+(1-s)\ty) (v-\tv)\cdot \nabla_y \rho_\vep\, ds
	\\
	&\qquad  \qquad - \int_0^1 D_yH(sv+ (1-s)\tv, sy+(1-s)\ty)(y-\ty)  \cdot \nabla_y \rho_\vep\, ds
	\\ & \geq \mu   \int_0^1 D_v H(sv+ (1-s)\tv, sy+(1-s)\ty) (y-\ty)\cdot \nabla_y \rho_\vep - (\mu+ \lambda) \ell_H\, (|\xi-\tilde \xi|  + |y-\ty|) 
	\\ & \geq   \mu\, \gamma \lambda \frac{ |y-\ty|^2}{\sqrt{\vep^2+ |y-\ty|^2}}- \mu^2\ell_H |y-\ty| - (\mu+ \lambda) \ell_H\, (|\xi-\tilde \xi|  + |y-\ty|)\,,
\end{align*}
where we have used $D_v H \geq\gamma\, I_d,$  {$||D_v H||, ||D_y H||\leq \ell_H$}, and  $|\nabla_y \rho_\vep - \lambda\frac{  (y-\ty) }{\sqrt{\vep^2+ |y-\ty|^2}}  | \leq \mu$.
Being 
$$
	\frac{ |y-\ty|^2}{\sqrt{\vep^2+ |y-\ty|^2}}\geq |y-\ty|_\vep \geq |y-\ty|- \vep\,,
$$
in view of the Lipschitz character of $B$  we get 
\begin{align*}
	& - (H(v,y)-H(\tv,\ty)) \cdot \nabla_y \rho_\vep  + ( B(v,y) -B(\tv,\ty) ) \cdot \nabla_v \rho_\vep
	\\ & \geq \mu\, \gamma \lambda \, |y-\ty| - \mu^2\ell_H |y-\ty| - (\mu+ \lambda) \ell_H\, (|\xi-\tilde \xi|  + |y-\ty|) - \mu\, \gamma\,  \lambda\, \vep   -     \ell_B( |v-\tv|+ |y-\ty|) 
	\\ & \geq ( \mu\, \gamma  \lambda- \mu^2\ell_H- (\mu+ \lambda) \ell_H- \ell_B (\mu+1) ) |y-\ty|- ((\mu+ \lambda) \ell_H+ \ell_B)  |\xi-\tilde \xi|- \mu\, \gamma\,  \lambda\, \vep\,.
\end{align*}
We recall that $\mu=  \frac{2\ell_H}{\gamma}\geq 2$, and then $ \mu\, \gamma  \lambda-  \lambda  \ell_H=\frac12 \mu\gamma\lambda$. Also using $ |\xi-\tilde \xi|\leq \sqrt \vep$, it holds
\begin{align*}
	& ( \mu\, \gamma  \lambda- \mu^2\ell_H- (\mu+ \lambda) \ell_H- \ell_B (\mu+1) ) |y-\ty|- ((\mu+ \lambda) \ell_H+ \ell_B)  |\xi-\tilde \xi|- \mu\, \gamma\,  \lambda\, \vep \\ & \geq 
	\mu\bigg( \frac12  \gamma \lambda  - (\mu +1) \ell_H  - 2 \ell_B  \bigg) |y-\ty|- C\sqrt\vep  = \frac14 \mu\gamma \lambda |y-\ty|- C\sqrt\vep\,,  
\end{align*}
by definition of $\lambda$. According to $ \mu\geq 2$ and $\lambda |y-\ty| \geq  \rho -   \sqrt\vep \geq \rho_\vep-   \sqrt\vep$, we conclude that 
$$
	- (H(v,y)-H(\tv,\ty)) \cdot \nabla_y \rho_\vep  + ( B(v,y) -B(\tv,\ty) ) \cdot \nabla_v \rho_\vep \geq \frac\gamma 2 \rho_\vep - C\sqrt\vep\,. 
$$
Inserting the above estimate  in \rife{zeta3}, we obtain
\begin{align*}
	   \left(K[\vfi( v,y)+ \vfi( \tv,\ty)]+ L\right)  \psi'(\rho_\vep) \left\{ \frac\gamma2 \rho_\vep - C \sqrt\vep \right\}  & \leq \psi(\rho_\vep) K\, 2k_0    + \omega \,L \psi(\rho_\vep) 
	   \\
	   & \leq \frac\gamma 4 \psi'(\rho_\vep)\rho_\vep\, L  + \left( 1+ \frac \gamma 4\right) \omega \,L \psi(\rho_\vep) 
\end{align*}
in view of \rife{quelk0}. Therefore, we infer  that 
\be\label{apr2}
	\psi'(\rho_\vep) \left\{  \frac\gamma4    \rho_\vep - C \sqrt\vep \right\} \leq \left( 1+ \frac \gamma 4\right)\omega\, \psi(\rho_\vep)\,.
\ee
Finally, we can choose $\omega$ so that
$$
	\left( 1+ \frac \gamma 4\right) \omega <    \frac\gamma8 \theta\,  e^{-C_2 R_1^\theta}   \,,
$$
in addition to \rife{omega1}. Since $\rho_\vep\leq R_1$, this yields
$$
	\left( 1+ \frac \gamma 4\right)\omega\, \psi(\rho_\vep)\leq  \frac\gamma8 \theta\,  e^{-C_2 \rho_\vep^\theta}\psi(\rho_\vep)\leq \frac\gamma8 \psi'(\rho_\vep)\rho_\vep\,,
$$
where we have exploited $\psi(r)\leq C_2 r^\theta$ and  \rife{psi'-times-r}.  
Then, by virtue of this  inequality and \rife{apr2}, we achieve 
\be\label{step2bis}
	\frac\gamma8    \rho_\vep \leq C  \sqrt\vep \,.
\ee
On the other hand, since $u(t )$ is Lipschitz continuous, with some constant $\ell_u^T$, uniformly   in $(0,T)$, by definition \rife{defmeps} and using \rife{rho_eps-bigger-sum-x-y} and \rife{step2bis}, we have 
\be\label{theend}
	{\mathbf m}\leq  {\mathbf m}_\vep \leq   u(t,v,y)- u(t, \tv, \ty )\leq \ell_u^T (\rho_\vep + C\vep)\leq \ell_u^T( 8C\gamma^{-1} \sqrt \vep + C\vep)\,.
\ee
Letting $\vep\to 0$, we obtain a contradiction because $\mathbf m>0$.

We have thus proved that \rife{claim} holds for $K$ given by \rife{choiK}, with a choice of $R_1, \omega,$ and $C_2$ only depending on $\ell_B, \ell_H, \gamma$ and the Lyapunov function $\vfi$. The constant $L$ is finally chosen in a way that  $K$, as defined by \rife{choiK}, satisfies \rife{Ku0}. This means taking $L=  n\, [u_0]_\theta$, for some $n\geq \frac{\omega_0 \vfi_{R_1}}{4( \omega + \tilde C\psi'(R_1)R_1) }$.
Letting $\eta\to 0$ in \rife{claim}, we obtain
$$
	u(t,v,y)- u(t, \tv, \ty ) \leq  e^{-\omega t}  \left(K[\vfi(v,y)+ \vfi(\tv,\ty)]+ L\right) \psi(\rho)  \leq e^{-\omega t}  (K+L) [\vfi(v,y)+ \vfi(\tv,\ty)] \psi(\rho)\,.
$$
Since $\psi(\rho)\leq (C_2 \rho^\theta)\wedge 1$ and $\rho$ is equivalent to the Euclidean distance, we deduce that $[{u(t)}]_\theta\leq C\, e^{-\omega t}  (K+L)$. Recalling that $K$ is proportional to $L$, where $ L=  n\, [u_0]_\theta$,  we conclude with the desired estimate. Hence, the case $\theta \in (0,1)$ is completed.
\vskip1em
Finally, we  explain how to get the result for the limiting case $\theta=1$. Note that we only need to refine the previous estimate when $\rho$ is small.  So, we only consider the case that $\rho<\de$, for a small $\de <1$ (to be fixed later). First, we modify the function $\psi$ used before. To this end, we fix any $\theta<1$, we take a constant $\beta$ (to be fixed later, depending on  $\de$) and define
$$
	\tilde \psi(r){:=} r - \frac\beta{1+\theta}r^{1+\theta}\,.
$$
We suppose, for now, that 
\be\label{betade1}
	\beta \de^\theta<\frac12
\ee
so that
\be\label{psi-lin}
 	\frac12 r\leq \tilde \psi(r)\leq r \,,\qquad \frac12 \leq \tilde \psi'(r) \leq 1  \qquad \forall r \leq \de\,.
\ee
We claim that, for a suitable choice of $\omega$ and $K$, it holds 
\be\label{claim2}
	u(t,v,y)- u(t, \tv, \ty ) \leq  e^{-\omega t}   K\, [\vfi(v,y)+ \vfi(\tv,\ty)]  \tilde \psi(\rho) + \frac{\eta}{T-t}\qquad \forall t\in (0,T),\,\,v,y,\tv,\ty\in \R^d\,:\, \rho<\de\,,
\ee 
for every $\eta$ sufficiently small. Here, $\vfi$ and $\rho$ are the same quantities used above, except that the parameters $\mu$ and $\lambda$ in \rife{rotated-norm} might be chosen differently (as precised later), yet only depending  on $\gamma, \ell_H,$ and $\ell_B$. Following the same approach as before, \rife{claim2} is proved arguing by contradiction, supposing that 
$$
	{\mathbf m}:= \sup_{\rho<\de}\,\,  \bigg\{ u(t,v,y)- u(t, \tv, \ty ) -   e^{-\omega t}  K[\vfi(v,y)+ \vfi(\tv,\ty)] \tilde  \psi(\rho) -  \frac{\eta}{(T-t)} \bigg\}    >0\,.
$$
If $\omega_\theta$ and $K_\theta$ denote the constants  given by \rife{ktheta}, here we choose $\omega<\omega_\theta$  and $K$ such that
$$
	K\geq  2\max\left(1, \de^{\theta-1}K_\theta\right)[u_0]_1\,.
$$
Then, { according to this and \rife{psi-lin}, when $t=0$, we have}
$$
	u_0(v,y)- u_0(\tv, \ty )   \leq  [u_0]_1[\vfi(v,y)+ \vfi(\tv,\ty)] \rho 	\leq  K[\vfi(v,y)+ \vfi(\tv,\ty)] \tilde  \psi(\rho) \qquad \forall \rho\leq \de \,,
$$
and, if $\rho=\de$, 
\begin{align*}
	u(t,v,y)- u(t, \tv, \ty ) & \leq   e^{-\omega_\theta t}  K_\theta [u_0]_\theta \, [\vfi(v,y)+ \vfi(\tv,\ty)] \de^\theta \leq e^{-\omega  t}  K_\theta \de^\theta [u_0]_1 \, [\vfi(v,y)+ \vfi(\tv,\ty)]  
	\\ & \leq e^{-\omega  t}  K {\frac\de2}  \, [\vfi(v,y)+ \vfi(\tv,\ty)]  \leq  e^{-\omega  t}  K [\vfi(v,y)+ \vfi(\tv,\ty)]\tilde \psi(\de)\,.
\end{align*}
Consequently,  if ${\mathbf m}$ is positive, then it can be attained neither at $t=0$ nor for $\rho=\de$. Excluding the case $t=T$ as before, it follows that ${\mathbf m}$ is a local maximum.  Then, we proceed exactly as in the first part of the proof. Specifically, we introduce the parameter $\vep$ and  the corresponding ${\mathbf m}_\vep$, which, for $\vep$ sufficiently small, is attained at some $t\in (0,T)$, and $ v,y,\tv,\ty$ with $\rho_\vep<\de$.  With the notation already used above, if $|\xi-\tilde \xi|>\sqrt \vep$, then we employ the coupling by reflection in the $\xi$-variable, and we obtain the analogue of \rife{nick}, namely,
\begin{align*}
	&  -2 \tilde \psi(\rho_\vep)k_0\,  K      -                     K[\vfi( v,y)+ \vfi( \tv,\ty)]  C\,   \tilde\psi'(\rho_\vep) \rho_\vep  \\
	& \qquad \leq    K[\vfi( v,y)+ \vfi( \tv,\ty)] 
	\bigg\{ 4     \,  \tilde\psi''(\rho_\vep)\, \bigg( \frac{|\xi-\tilde \xi|}{|\xi-\tilde \xi|_\vep + \vep}\bigg)^2 
	+ 4     \,\frac{ \tilde\psi'(\rho_\vep)}{|\xi-\tilde \xi|_\vep + \vep}\, \frac{\vep^2}{(|\xi-\tilde \xi|_\vep + \vep)^2}\bigg\}
	\\
	& \qquad \leq    K[\vfi( v,y)+ \vfi( \tv,\ty)] 
	\bigg\{ 2     \,  \tilde\psi''(\rho_\vep) 
	+ 4     \sqrt \vep \tilde \psi'(\rho_\vep) \bigg\}\,,
\end{align*}
by virtue of \rife{rooteps} and the concavity of $\tilde \psi$ again. We recall that  $C$ only depends on $\gamma, \ell_H,\ell_B,$ and $\ell_\vfi$.   
Using \rife{psi-lin}, we deduce
$$
	0 \leq    K[\vfi( v,y)+ \vfi( \tv,\ty)] 
	\bigg\{ 2     \,  \tilde\psi''(\rho_\vep) 
	+ 4     \sqrt \vep +(2k_0  + C ) \rho_\vep\bigg\}\,,
$$
whence, computing $\tilde\psi''$ and being $|\xi-\tilde \xi|>\sqrt \vep$, 
$$
	2\theta \beta \rho_\vep^{\theta-1}\leq 4\sqrt \vep + (2k_0  + C ) \rho_\vep\leq \tilde C \rho_\vep\,,
$$
for some $\tilde C$ only depending on $\gamma, \ell_H,\ell_B,$ and $\vfi$. We now choose $\beta=  \frac{\tilde C}{\theta} \de^{2-\theta} $. In this way, the previous inequality cannot occur for $\rho_\vep \leq \de$. Hence,  the case $|\xi-\tilde \xi|>\sqrt \vep$ is not possible. Observe that, with this choice of $\beta$,  \rife{betade1} is satisfied as long as $\de$ is sufficiently small. It remains to deal with $|\xi-\tilde \xi|\leq \sqrt \vep$. We reason as in Case 3 above. In particular, we estimate
\begin{align*}
	& - (H(v,y)-H(\tv,\ty)) \cdot \nabla_y \rho_\vep  + ( B(v,y) -B(\tv,\ty) ) \cdot \nabla_v \rho_\vep
	\\ & \geq ( \mu\, \gamma  \lambda- \mu^2\ell_H- (\mu+ \lambda) \ell_H- \ell_B (\mu+1) ) |y-\ty|- ((\mu+ \lambda) \ell_H+ \ell_B)  |\xi-\tilde \xi|- \mu\, \gamma\,  \lambda\, \vep
	\\ & \geq ( \mu\, \gamma  \lambda- \mu^2\ell_H- (\mu+ \lambda) \ell_H- \ell_B (\mu+1) ) |y-\ty|- C\sqrt \vep \,,
\end{align*}
which yields, from $ \rho_\vep {+\lambda\sqrt{\vep}}\geq \lambda |y-\ty|\geq   \rho_\vep - \sqrt \vep$,   
\begin{align*}
	& - (H(v,y)-H(\tv,\ty)) \cdot \nabla_y \rho_\vep  + ( B(v,y) -B(\tv,\ty) ) \cdot \nabla_v \rho_\vep
	\\ & \geq ( \mu\, \gamma- \ell_H) \rho_\vep     -  \frac1\lambda (\mu^2\ell_H+ \mu\ell_H + \ell_B (\mu+1) )\rho_\vep   - C\sqrt \vep \,,
\end{align*}
for a possibly different $C$, potentially depending on $\mu$ and $\lambda$. The analogue of \rife{zeta3}, in this case, leads to 
\begin{align*}
	( \mu\, \gamma- \ell_H) \rho_\vep     -  \frac1\lambda (\mu^2\ell_H+ \mu\ell_H + \ell_B (\mu+1) )\rho_\vep  
	&  \leq 2 \frac{\psi(\rho_\vep)}{\psi'(\rho_\vep)} k_0\,       + C\sqrt \vep  
	\\
	&   \leq 4 \rho_\vep k_0\,        + C\sqrt \vep \,,
\end{align*}
also using \rife{psi-lin}. Now, we fix the values of $\mu$ and $\lambda$.  First, we  choose $\mu = \frac{4k_0+ 2\ell_H}\gamma$ so that the previous inequality yields 
$$
	\ell_H  \rho_\vep     -  \frac1\lambda (\mu^2\ell_H+ \mu\ell_H + \ell_B (\mu+1) )\rho_\vep \leq C\sqrt \vep \,.
$$
Next, we choose $\lambda $ large to reach that  $\rho_\vep \leq \hat C \sqrt \vep$, for some different $\hat C$. From here, we  obtain \rife{theend} and conclude the contradiction argument as before.
\end{proof} 
\qed

{ \begin{remark}\label{extensions} The approach developed in the above proof can be readily adapted to slightly different settings. For instance, if the state variable $y$ belongs to a compact space, like the $d$-dimensional torus, then it can be enough to use a Lyapunov function $\vfi$ only depending on the $v$-variable; in that case, once the possibility that $|v-y|$ is large is excluded by the confining role of $\vfi$, the rest of the proof is reduced to a compact framework for the whole variable $(v,y)$.

A less trivial extension, though valuable, can be obtained if the variable $y$ belongs to a space with a different dimension, say $y\in \R^k$, $k<d$. There, we need to change the scalar $\mu$ to a suitable matrix $M\in {\mathcal  M}_{k,n}$ in the definition of the rotated norm, and the nondegeneracy condition on $D_vH$ to a full-rank condition.
\end{remark}
}

\subsection{Decay of Wasserstein distances for  the  kinetic Fokker-Planck equation}

An immediate corollary is the following decay result {in Wasserstein distances} for solutions of \rife{FP}. Let us introduce the following weighted  distance for $\mu, \nu \in { \cM_\vfi(\R^{2d})\cap \cP_1({\R^{2d}})}$: 
$$
	d_{1,\vfi}(\mu, \nu):= \sup\left\{ \into \zeta\, d(\mu-\nu)\,,\quad \zeta\in L^\infty(\vfi^{-1})\,:\, \, [\zeta]_1\leq 1\right\}\,,
$$
{ where $[\cdot]_1$ is defined in \rife{def-seminorm}. Then, for every $\zeta\in L^\infty(\vfi^{-1})$, one can estimate
\be\label{tauto}
\into \zeta\, d(\mu-\nu) \leq [\zeta]_1\, d_{1,\vfi}(\mu, \nu)\,.
\ee}
Note that if $\vfi$ has at least linear growth at infinity, then we have  
$$
\frac{|u(\chi)-u(\tilde \chi)|}{(\vfi(\chi)+ \vfi(\tilde \chi)) (|\chi-\tilde \chi|\wedge 1)}\leq C\, \frac{|u(\chi)-u(\tilde \chi)|}{|\chi-\tilde \chi |}\,, \quad \chi= (v,y)\,,\,\,\tilde \chi=(\tilde v, \tilde y)\,,
$$
for some $C>0$ independent of $u$.  Hence, we get 
\be\label{wass-dist}
	\cW_1(\mu, \nu)\leq C\, d_{1,\vfi}(\mu, \nu)\,,
\ee
where $\cW_1$ is the standard $1$-Wasserstein distance, {recall \rife{Was-dist}}.  In particular,  in this case, the decay of $d_{1,\vfi}$ implies the decay of $\cW_1$.

\begin{theorem}\label{decay-Wass} Under the assumptions of Theorem \ref{LTO}, let $\mu_1,\mu_2$ be two solutions of \rife{FP} corresponding to initial data $\mu_{01}, \mu_{02}\in \cP_1(\R^{2d})$.  Then, there exist $K, \omega$ such that it holds
\be\label{dirate}
d_{1,\vfi}({\mu_1(t)}, {\mu_2(t )}) \leq K e^{-\omega t } d_{1,\vfi}(\mu_{01}, \mu_{02})\,.
\ee
In particular, if $\vfi$ grows at least linearly at infinity, then we have
$$
	\cW_1({\mu_1(t )}, {\mu_2(t )}) \leq C\, K\,  e^{-\omega t } d_{1,\vfi}(\mu_{01}, \mu_{02})
$$
for some $C>0$.
\end{theorem}

\begin{proof} Let $\zeta\in  L^\infty(\vfi^{-1})$ be such that $[\zeta]_1\leq 1$. This means that $\zeta$ is a  locally Lipschitz function and { (see Corollary \ref{visco-unbounded})} there exists a viscosity solution $u$ of the backward problem
$$
	\begin{cases}
	-\partial_t u+ \cL[u]     = 0 &  \text{ in } ( 0,t ) \times \R^{2d}\,, \\ 
	u(t) = \zeta   &  \text{ in }   \R^{2d}\,.
	\end{cases} 
$$
By \rife{dual-form} applied to both $\mu_1$ and $\mu_2$, it holds
$$
	\into \zeta\, d({\mu_1(t )- \mu_2(t )})= \into u(0)\, d(\mu_{01}- \mu_{02}) \leq Ke^{-\omega t} d_{1,\vfi}(\mu_{01}, \mu_{02})\,,
$$
where we used \rife{tauto}  and  \rife{ktheta}  (with $\theta=1$) to estimate  $[u(0)]_1$. This   concludes with \rife{dirate}.  The last statement immediately follows from \rife{wass-dist}.
\end{proof}
\qed
\vskip0.5em

We conclude by pointing out that Theorem \ref{decay-Wass} provides the existence (and uniqueness) of a stationary measure $\bar m$ solving \rife{FP} and attracting all other solutions in long time. As detailed in \cite[Section 5]{Po}, the exponential decay rate enables building such a solution, whose uniqueness and exponential stability follow from the same contractivity estimate.

\section{Hypocoercivity and $\mathbf {L^\infty}$-Lipschitz smoothing}\label{sect-hypo}

In this section, we prove short-time smoothing  estimates for the problem \rife{pb}.  To work in a slightly more general framework, we consider the operator
$$
	\cL (u):=  - {\rm tr}(A(v)D^2_{vv} u) - \kappa\,  \Delta_y u  - H(v,y) \cdot \nabla_y u +  B(v,y) \cdot \nabla_v u\,,
$$
where $A(v), H(v,y),$ and $B(v,y)$ are $C^1$-functions, globally Lipschitz continuous, and $\kappa \geq 0$. 
In particular, we assume that $A(v)$ is a bounded and elliptic matrix, i.e., 
\be\label{A1}
	\exists \,\, \alpha_0>0\,: \quad  \alpha_0 I_d \leq A(v) \leq \frac1{\alpha_0} I_d\,,
\ee
and  also Lipschitz continuous:
\be\label{Alip}
	| A(v)-A(\tv )| \leq \ell_A \, |v-\tilde v|\qquad \forall x,\tv \in \R^d\,.
\ee
We  start with a computational lemma. 

\begin{lemma}\label{calc-hyp} 
If $u$ is a smooth solution of 
$$
	\partial_t u+ \cL(u) =0\,,
$$
then  we have:
\begin{itemize}
\item[(i)] 
$$
	\frac12 (\partial_t + \cL) \left( u^2\right)= - A(v) \nabla_v u \cdot \nabla_v u - \kappa\, |\nabla_y u|^2\,;
$$

\item[(ii)] 
\begin{align*}
	\frac12 (\partial_t + \cL) \left( |\nabla_v u|^2\right) & = - \sum_k A(v) \nabla_v u_{v_k}  \cdot  \nabla_v u_{v_k} + \sum_{ijk} \, (a_{ij})_{v_k}u_{v_k} u_{v_iv_j} 
	\\
	& \, + 
	D_v H\nabla_y u  \cdot \nabla_v u - D_v B \nabla_vu \cdot \nabla_v u - \kappa\, \sum_{ij} u_{v_iy_j}^2\,;
\end{align*}
\item[(iii)] 
\begin{align*}
	(\partial_t + \cL) & \left(  \nabla_v u\cdot \nabla_y u\right)   =\sum_{ijk} \, ( (a_{ij})_{v_k}u_{y_k} u_{v_iv_j}- 2 a_{ij} u_{v_iv_k} u_{y_kv_j} )-2\kappa\,  \sum_{ij} u_{v_iy_j} u_{y_iy_j}
	\\
	& \, + D_v H \nabla_y u\cdot \nabla_y u  +D_y H \nabla_y u \cdot  \nabla_v u - D_v B\nabla_v u\cdot \nabla_y u - D_y B \nabla_v u\cdot \nabla_v u\,;  
\end{align*}
\item[(iv)]
$$
	\frac12 (\partial_t + \cL)  \left( |\nabla_y u|^2\right)  =- \sum_k A(v) \nabla_v u_{y_k}  \cdot \nabla_v u_{y_k} + D_y H \nabla_y u\cdot \nabla_y u
	- D_y B \nabla_v u\cdot  \nabla_y u - \kappa\,  \sum_{ij} u_{y_iy_j}^2 \,.
$$
\end{itemize}
\end{lemma}

\begin{proof}  Equality (i)   immediately follows from the chain rule and the equation  $\partial_t u+ \cL(u) =0$. To show (ii), we denote the second-order terms in the operator $\cL$ by $\cL_0 (u):=  - {\rm tr}(A(v)D^2_{vv} u)- \kappa\, \Delta_y u$. Then, we have 
\begin{align*}
	& \frac12 (\partial_t + \cL_0)\left( |\nabla_v u|^2\right)   = \nabla_v  \left( (\partial_t + \cL_0)(u)\right) \cdot \nabla_v u + \sum_{ijk} \, \left((a_{ij})_{v_k}u_{v_k} u_{v_iv_j}  - a_{ij} u_{v_kv_i}u_{v_kv_j} \right) - \kappa\, \sum_{ij} u_{v_iy_j}^2
	\\ & \quad = \nabla_v \left(  H(v,y) \cdot \nabla_y u -  B(v,y) \cdot \nabla_v u\right) \cdot \nabla_v u + \sum_{ijk} \, \left((a_{ij})_{v_k}u_{v_k} u_{v_iv_j} -  a_{ij} u_{v_kv_i}u_{v_kv_j}\right) - \kappa\, \sum_{ij} u_{v_iy_j}^2\,,
\end{align*}
in view of $\partial_t u+ \cL(u) =0$ again. 
Developing the drift terms, we reach (ii).

Similarly,  we have
$$
	\frac12 (\partial_t + \cL_0)\left( |\nabla_y u|^2\right)   = \nabla_y  \left( (\partial_t + \cL_0)(u)\right) \cdot \nabla_y u  - \sum_{ijk} \, a_{ij} u_{y_kv_i}u_{y_kv_j}  - \kappa  \sum_{ij} u_{y_iy_j}^2\,,
$$
and  we obtain (iv) using once more the equation for $u$.

Finally, in order to   get (iii), we compute
\begin{align*}
	& (\partial_t + \cL_0)   \left( \nabla_v u\cdot \nabla_y u\right)   =\nabla_v  \left( (\partial_t + \cL_0)(u)\right) \cdot \nabla_y u   +  \nabla_y  \left( (\partial_t + \cL_0)(u)\right) \cdot \nabla_v u  \\ & \qquad+ \sum_{ijk} \,  \left( (a_{ij})_{v_k}u_{y_k} u_{v_iv_j} -2 a_{ij} u_{v_iv_k}u_{v_jy_k} \right) - 2\kappa\,  \sum_{ij} u_{v_iy_j} u_{y_iy_j}
	\\ &=  \nabla_v \left(  H(v,y) \cdot \nabla_y u -  B(v,y) \cdot \nabla_v u\right)\cdot  \nabla_y u + \nabla_y \left(  H(v,y) \cdot \nabla_y u -  B(v,y)\cdot\nabla_v u\right) \cdot \nabla_v u  \\ & \qquad+ \sum_{ijk} \,  \left( (a_{ij})_{v_k}u_{y_k} u_{v_iv_j} -2 a_{ij} u_{v_iv_k}u_{v_jy_k} \right)- 2\kappa\,  \sum_{ij} u_{v_iy_j} u_{y_iy_j}
	\\  & = H(v,y)\cdot \nabla_y (\nabla_v u \cdot \nabla_y u) - B(v,y) \cdot \nabla_v (\nabla_v u \cdot \nabla_y u) +   D_v H \nabla_y u \cdot  \nabla_y u+D_y H \nabla_y u \cdot  \nabla_v u \\
	& \qquad  - D_v B \nabla_v u\cdot \nabla_y u - D_y B \nabla_v u\cdot \nabla_v u    + \sum_{ijk} \,  \left( (a_{ij})_{v_k}u_{y_k} u_{v_iv_j} -2 a_{ij} u_{v_iv_k}u_{v_jy_k} \right) - 2\kappa\,  \sum_{ij} u_{v_iy_j} u_{y_iy_j}\,{.}  
\end{align*}
\end{proof}
\qed

\vskip1em

Exploiting the above lemma, we deduce a short-time regularization with an $L^\infty$-version of  hypocoercivity estimates.

\begin{proposition}\label{hypocoer} Let $u\in L^\infty((0,T); W^{1,\infty}({\R^{2d}}))$  be a smooth solution of the equation
$$
	\begin{cases}
		\partial_t u - {\rm tr}(A(v)D^2_{vv} u) - \kappa\,\Delta_y u -  H(v,y)  \cdot \nabla_y u +  B(v,y)\cdot\nabla_v u=0\\
		u(0)=u_0\,,
	\end{cases}
$$
where $\kappa\geq 0$, $A(v)$ satisfies \rife{A1}-\rife{Alip},  $B(v,y)$ fulfils  \rife{Blip}, and $H(v,y)$ satisfies \rife{H1}-\rife{H2}.

\vskip0.6em
\noindent Then, there exists a constant $C>0$, depending on $\alpha_0, \ell_A,  \ell_B, \ell_H,$ and $\gamma$  such that
\be\label{hypo}
	\, \|\nabla_v u{(t)}\|_\infty + t\,  \|\nabla_y u{(t)}\|_\infty \leq  C\, \frac{\|u_0\|_\infty}{\sqrt t } \qquad \forall t\in (0,1]\,.
\ee
\end{proposition}

\begin{proof} We introduce the function
\be\label{Vilfunction}
W: = \frac12 \left( \lambda |u |^2 + t |\nabla_v u|^2 - 2\vep t^2 \nabla_v u \cdot \nabla_y u + \de t^3 |\nabla_y u|^2\right),
\ee
where $\lambda$ is a large constant, to be fixed later, and $\vep, \de$ satisfy (at least) the condition
\be\label{vepde1}
\vep^2 < \de\,,
\ee
so that $W>0$. The goal is to show that $W $ is a subsolution (in short time).

From  Lemma \ref{calc-hyp}-(ii), using \rife{Blip} and the ellipticity of $A(v)$, see \rife{A1}, we have
\be\label{ber1}
	\begin{split}
		\frac12 (\partial_t + \cL) \left( t \, |\nabla_v u|^2\right)  & \leq \bigg(\frac12+ t\, \ell_B\bigg) |\nabla_v u|^2  -  t\, \alpha_0 \sum_{ik}   u_{v_iv_k}^2  + t\, \sum_{ijk} \, (a_{ij})_{v_k}u_{v_k} u_{v_iv_j}  \\ & \qquad +  t\, D_v H\nabla_y u  \cdot \nabla_v u - t\kappa\, \sum_{ij} u_{v_iy_j}^2  \\
		& \leq \bigg(\frac12+  t\, (\ell_B+ C_A ) \bigg) |\nabla_v u|^2  -  t\, \frac{\alpha_0}2 \sum_{ik}   u_{v_iv_k}^2 +  t\, \ell_H\,  |\nabla_v u| \, |\nabla_y u |- t\kappa\, \sum_{ij} u_{v_iy_j}^2\,,
	\end{split}
\ee
thanks to  \rife{Alip}, \rife{H2}, and Young's inequality. Hereafter, by $C_A$, we denote possibly different constants  depending on the $W^{1,\infty}$ bound of $A$, namely, on $\alpha_0$ and $\ell_A$.\\
Similarly,  according to Lemma \ref{calc-hyp}-(iii),  by virtue of  \rife{A1}--\rife{Alip} and Young's inequality again, it holds
\begin{align*}
	(\partial_t& + \cL)  \left(- \vep\, t^2\,  \nabla_v u\cdot \nabla_y u\right)   \leq   t\, \frac{\alpha_0}4  \sum_{ij} \, u_{v_iv_j}^2 +  {C_A}\, \vep^2\, t^3   |\nabla_y u|^2 +  {C_A}\,\vep^2\, t^3  \sum_k |\nabla_v u_{y_k}|^2  \\
	& - \vep \, t^2\,  D_v H \nabla_y u\cdot \nabla_y u  - \vep \, t^2\, D_y H \nabla_y u \cdot  \nabla_v u +  \vep\, t^2 (D_v B\nabla_vu\cdot \nabla_y u + D_y B \nabla_vu\cdot \nabla_v u)  
	\\ & 
	\qquad 
	+ 2\kappa\, \vep\,t^2 \sum_{ij} u_{v_iy_j} u_{y_iy_j} -2\vep\, t\, \nabla_v u\cdot \nabla_y u\,.
\end{align*}
The terms with $H$ and $B$ can be estimated from \rife{H1}--\rife{H2} and \rife{Blip}. So, we get
\begin{align*}
		(\partial_t + \cL)&\left(- \vep\, t^2\,  \nabla_v u\cdot \nabla_y u\right)  \leq  t\, \frac{\alpha_0}4  \sum_{ij} \, u_{v_iv_j}^2 +  {C_A}\, \vep^2\, t^3   |\nabla_y u|^2 +  {C_A}\,\vep^2\, t^3  \sum_k |\nabla_v u_{y_k} |^2  
		\\ & \qquad - \vep \frac\gamma2 \, t^2 |\nabla_y u|^2 + C_{B,H}\,  \vep\, t^2  |  \nabla_v u|^2+ 2\kappa\, \vep\,t^2 \sum_{ij} u_{v_iy_j} u_{y_iy_j} -2\vep\, t\, \nabla_v u\cdot \nabla_y u\,,
\end{align*}
where $C_{B,H}$ only depends on the constants $\ell_B, \ell_H, \gamma$.\\
Finally, in view of Lemma \ref{calc-hyp}-(iv), we also have
\begin{align*}
		\frac12 (\partial_t + \cL) \left(\de\, t^3\,  |\nabla_y u|^2\right)&  \leq  -\de\, t^3\alpha_0  \sum_k |\nabla_v u_{y_k} |^2  + \de \, t^3  \left( \ell_B| \nabla_v u |\, |\nabla_y u|+ \ell_H |\nabla_y u|^2\right) 
		\\ & \qquad 
		- \kappa\, \de\, t^3\, \sum_{ij} u_{y_iy_j}^2+ \frac32 \de\, t^2 \,  |\nabla_y u|^2\,.
\end{align*}

We now add up the previous inequalities, gathering together the similar terms.  Hence, we obtain 
\begin{align*}
	(\partial_t + \cL) & \left( W- \frac \lambda2  u^2 \right)   \leq  \bigg(\frac12+  t\, (\ell_B+ C_A ) \bigg) |\nabla_v u|^2  -  t\, \frac{\alpha_0}4 \sum_{ik}   u_{v_iv_k}^2 +  t\, \ell_H\,  |\nabla_v u| \, |\nabla_y u | 
	\\ & -   t^2 |\nabla_y u|^2 \left( \vep \frac\gamma2- \frac32 \de -  {C_A}\, \vep^2\, t\right)  -  t^3\left( \de\alpha_0 - {C_A}\vep^2\right) \sum_k |\nabla_v u_{y_k} |^2  
	   -2\vep\, t\, \nabla_v u\cdot \nabla_y u
	\\ & \quad + C_{B,H}\,  \vep\, t^2  |  \nabla_v u|^2+ \de \, t^3  \left( \ell_B| \nabla_v u |\, |\nabla_y u|+ \ell_H |\nabla_y u|^2\right) 
	\\ & \qquad 
	+ \kappa\,{t}\, \left( 2\vep\,t \sum_{ij} u_{v_iy_j} u_{y_iy_j}-  \de\, t^2\, \sum_{ij} u_{y_iy_j}^2 -  \sum_{ij} u_{v_iy_j}^2\right)\,.
\end{align*}
Recall that we assume  $t\leq 1$ and $\vep, \de\leq 1$. Moreover, thanks to \rife{vepde1}, the last term is negative and can be dropped.  Then, applying Young's inequality, we infer that there exists $K>0$ (only depending on $\alpha_0, \ell_A, \ell_B, \ell_H,$ and $\gamma$) such that
 \begin{align*}
	(\partial_t + \cL) & \left( W- \frac \lambda2  u^2 \right)   \leq  K\,  |\nabla_v u|^2  
+  (\ell_H+2\vep) t\, |\nabla_v u| \, |\nabla_y u | 
	\\ & \quad - t^2 \left( \vep \frac\gamma2- \frac32 \de -  {C_A}\, \vep^2\, t-   \ell_B\de^2 t- \ell_H \de\, t\right)    |\nabla_y u|^2 -  t^3\left( \de\alpha_0 - {C_A}\vep^2\right)  \sum_k |\nabla_v u_{y_k} |^2 \,,
\end{align*}
which yields, { because $(\ell_H+2\vep) t\, |\nabla_v u| \, |\nabla_y u | \leq \frac{(\ell_H+2\vep)^2}{ \vep\gamma} |\nabla_v u|^2+ t^2 \vep\,\frac\gamma 4 |\nabla_y u|^2$,
\be\label{ber4}
 \begin{split}
	(\partial_t + \cL) & \left( W- \frac \lambda2  u^2 \right)   \leq  \left( K + \frac{(\ell_H+2\vep)^2}{ \vep\gamma} \right)  |\nabla_v u|^2   
	\\ & \quad - t^2 \left( \vep \frac\gamma4- { \frac{(3+ 2\ell_H)}2} \de -  {C_A}\, \vep^2 -   \ell_B\de^2 \right)    |\nabla_y u|^2 -  t^3\left( \de\alpha_0 - {C_A}\vep^2\right)  \sum_k |\nabla_v u_{y_k}|^2\,,  
\end{split}
\ee
where we have also used  $t\leq 1$ inside the second term on the right.}
We remark that the last two terms are negative as long as $\de, \vep$ are sufficiently small with $2(3+2\ell_H)\de<\gamma\vep$ { and $\de\alpha_0>C_A\vep^2$. }  We fix, e.g.,  $\de= \frac{\gamma\vep}{4(3+2\ell_H)}$ and take  $\vep$ small enough, to ensure  that \rife{vepde1} is also satisfied. Then, { being $\de=O(\vep)$, } it holds
\begin{align*}
{ - t^2 \left( \vep \frac\gamma4- { \frac{(3+ 2\ell_H)}2} \de -  {C_A}\, \vep^2 -   \ell_B\de^2 \right)  \leq  - t^2 \left( \vep \frac\gamma8  + O(\vep^2)  \right) }\,,
\end{align*}
{ and the right-hand side is negative choosing $\vep$ sufficiently small, {only depending on $C_A$ and $\ell_B$}.}  We can argue similarly for the last term in \rife{ber4}, { which is negative by taking $\vep$ small enough. Therefore,  we achieve}
$$
	(\partial_t + \cL)   \left( W- \frac \lambda2  u^2 \right)   \leq  \left( K + \frac{{(\ell_H+2\vep)^2}}{ \vep\gamma} \right)  |\nabla_v u|^2\,,   
$$
and, recalling Lemma \ref{calc-hyp}-(i), we reach
$$
	(\partial_t + \cL)  W \leq -  \left(  \lambda\, \alpha_0 - \left(K + \frac{{(\ell_H+2\vep)^2}}{ \vep\gamma} \right) \right)  |\nabla_v u|^2\qquad \forall t\in (0,1]\,.
$$
Here, we choose $\lambda$ sufficiently large so that $(\partial_t + \cL)  W\leq 0$ for $t\leq 1$. Now, $W$ is a bounded  subsolution (since $u\in L^\infty((0,T); W^{1,\infty}(\R^d))$), so we can apply the maximum principle to get 
$$
	\sup_{t\in (0,1]} W{(t )} \leq \| W{(0 )}\|_\infty \leq \frac\lambda 2 \|u_0\|_\infty^2\,.
$$
In view of \rife{vepde1}, we have $W\geq c_\vep \left(t|\nabla_v u{(t )} |^2+ t^3|\nabla_y u{(t )}|^2\right)$, from which estimate \rife{hypo} follows.
\end{proof}
\qed

\begin{remark} { At different steps, one may use  some dissipativity property of $D_vB$, $D_yB$, $D_yH$ (which enter in (ii)-(iv) of Lemma \ref{calc-hyp}) in an attempt to improve the above estimates. However, here we are just exploiting the smoothing effect in {\it short-time}, so this seems unnecessary at this stage. In parallel, it could be interesting to understand under what conditions on the drift terms an estimate such as \rife{hypo} can hold globally in time.}
\end{remark}

We now observe that a weighted version of Proposition \ref{hypocoer} can be obtained similarly.

\begin{proposition}\label{hypocoer2} Let $u\in L^\infty((0,T); W^{1,\infty}({\R^{2d}}))$ be a smooth solution of the equation
$$
	\begin{cases}
		\partial_t u - {\rm tr}(A(v)D^2_{vv} u) - \kappa\,\Delta_y u -  H(v,y)  \cdot \nabla_y u +  B(v,y) \cdot \nabla_v u=0\\
		u(0 )=u_0\,, 
	\end{cases}
$$
where $A(v)$ satisfies \rife{A1}-\rife{Alip}, $B(v,y)$ fulfils \rife{Blip}, and $H(v,y)$ satisfies \rife{H1}-\rife{H2}, and $\kappa\geq 0$. { Assume there exists some function 
$\vfi\in C^{1,2}([0,T)\times \R^{2d})$  such that $\vfi\geq 1$,  $\vfi(t,v,y)\mathop{\longrightarrow}\limits^{|(v,y)|\to \infty} \infty$ uniformly for $t\in (0,T)$, and satisfies}  
\be\label{vfi2}
 (\partial_t + \cL)(\vfi^2)\geq - k_0\,,\qquad t\in (0,1),\,\, v,y\in \R^d\,,
\ee
for some $k_0>0$. Then, { we have $\|{u(t)}\|_{L^\infty(\vfi^{-1}) }\leq e^{\frac {k_0}2 t} \|u_0\|_{L^\infty(\vfi_0^{-1}) }$ for  any $t\in (0,1]$, and there exists  $C>0$, depending on $\alpha_0, \ell_A, \ell_B, \ell_H, \gamma,$ and $k_0$, such that
\be\label{hypo2}
\, |\nabla_v u(t, v,y)|  + t\,  |\nabla_y u(t,v,y)|  \leq  C\, {\vfi(t,v,y)}\, \frac{ \|u_0\|_{L^\infty(\vfi_0^{-1}) }}{\sqrt t } \qquad \forall t\in (0,1],\,\, v,y\in \R^d\,,
\ee
where $\vfi_0(\cdot)=\vfi(0, \cdot)$.}
\end{proposition}

{ \begin{remark} In this result, the Lyapunov function $\vfi$ may also depend on time, as we are only concerned with short-time estimates here. Indeed, given the linear growth of the drift fields $H$ and $B$, the function $\vfi= e^{\lambda t} ( \langle v\rangle^k+ M)$, with $\langle x\rangle= \sqrt{1+|v|^2}$, {satisfies \rife{vfi2} and the above conditions}, provided that $0<k\leq 2$ and $\lambda,M$ are chosen sufficiently large.
\end{remark}}

\begin{proof} Take $L:= \|u_0\|_{L^\infty(\vfi_0^{-1})}$, which means that $L= \sup \left(\frac {|u_0|}{\vfi_0}\right)$. We claim that 
$$
	e^{-k_0 t}W- \frac\lambda 2 L^2\vfi^2\leq 0\qquad \text{for } t\in (0,1]\,,
$$
where $W$ is the function \rife{Vilfunction} built in Proposition \ref{hypocoer}.
Of course, the above inequality is true at $t=0$ by definition of $L$ and $W$. Next, we point out that
\begin{align*}
	(\partial_t + \cL)\bigg(e^{-k_0 t}W- \frac\lambda 2L^2\vfi^2\bigg) & = - k_0 e^{-k_0 t}W+ e^{-k_0 t}(\partial_t + \cL)W-  \frac\lambda 2 L^2 (\partial_t + \cL)\vfi^2
	\\
	& \leq - k_0 e^{-k_0 t}W+  \frac\lambda 2 L^2k_0\,,
\end{align*}
where we have used the subsolution property of $W$ and {\rife{vfi2}}. If it happens that $e^{-k_0 t}W- \frac\lambda 2L^2\vfi^2 $ has a positive maximum, then we have $e^{-k_0 t} W>  \frac\lambda 2 L^2 \vfi^2 \geq \frac\lambda 2L^2$ at {a maximum} point, which yields  
$$
	(\partial_t + \cL)\bigg(e^{-k_0 t}W- \frac\lambda 2L^2\vfi^2\bigg) \leq - k_0 e^{-k_0 t}W+  \frac\lambda 2 L^2k_0  <0\,
$$
there. This violates the maximum principle, showing that $e^{-k_0 t}W- \frac\lambda 2 L^2\vfi^2 \leq 0$ for $t\in (0,1]$. Hence, $W\leq e^{k_0 t} \frac\lambda 2 L^2\, \vfi^2$, which gives the { desired conclusions  according to the definition of $W$. } 
\end{proof}
\qed

{ 
To conclude this section, we note that the above estimate allows us to extend the well-posedness result of Theorem \ref{teo-ex} to possibly unbounded initial data $u_0$.

\begin{corollary}\label{visco-unbounded} Assume that \rife{Blip} and \rife{H1}-\rife{H2} hold, and $\kappa\geq 0, T>0$.  Let $\vfi(t,v,y)$ be a function satisfying the assumptions of Proposition \ref{hypocoer2}. Given $u_0 \in C(\R^{2d})\cap L^\infty(\vfi_0^{-1})$, there exists a unique  viscosity solution $u\in C( [0,T]\times \R^{2d})\cap L^\infty(\vfi^{-1})$  of  the problem \rife{pb}. Moreover, for every $t>0$, {$u(t)$} is locally Lipschitz continuous and fulfils estimate \rife{hypo2}.
\end{corollary}

\begin{proof} The existence follows by approximation. Precisely, if $u_0{\coloneqq }\xi\, \vfi_0$, for some bounded $\xi$, then we consider the viscosity solutions $u_n$ corresponding to the initial data $u_{0n}= \xi (\vfi_0\wedge n)$. By Proposition \ref{hypocoer2},  we have 
$$
\|{u_n(t)}\|_{L^\infty(\vfi^{-1}) }\leq C_T \|u_{0n}\|_{L^\infty(\vfi_0^{-1}) } \leq C_T \|u_{0}\|_{L^\infty(\vfi_0^{-1}) }\,,
$$
and then by \rife{hypo2} we deduce uniform bounds for  $\sqrt t \, \|\nabla_v {u_n(t)}\|_{L^\infty(\vfi^{-1}) }$ and $t^{\frac32} \, \|\nabla_y {u_n(t)}\|_{L^\infty(\vfi^{-1}) }$. Henceforth, one proceeds as in the proof of Proposition \ref{viscos}  to deduce the  compactness of $u_n$ (in local uniform topology) and the convergence of a subsequence towards a viscosity solution $u$.  The uniqueness follows as in Proposition \ref{viscos}, replacing $|v|^2$ with $\vfi^2$ and using the supersolution property \rife{vfi2}. 
\end{proof}
\qed
}

\section{Improved long-time decay}\label{final-decay}

Thanks to Proposition \ref{hypocoer2},  we can now improve the long-time decay estimate of Theorem \ref{LTO}.

\begin{theorem}\label{LTO2}  Assume that \rife{Blip} and \rife{H1}-\rife{H2} hold, and $\kappa\geq 0$. { Let $\vfi(v,y)$ be a Lyapunov function (cf. Definition \ref{def-Lyap}) fulfilling \rife{liplocvfi} and \rife{vfi2}. In addition, we suppose}
\be\label{vfi3} 
	\int_0^1 \vfi(sv+(1-s)\tv, sy+(1-s)\ty) \, ds \leq c\, [\vfi(v,y)+ \vfi(\tv, \ty)] \qquad \forall v,y,\tv,\ty \in \R^d\,,
\ee
for some $c>0$. Then, there exist $\omega, K>0$ such that the viscosity solution of \rife{pb} satisfies
\be\label{goal}
| u(t,v,y ) - u(t, \tv, \ty) | \leq K\, e^{-\omega t}  \|u_0\|_{L^\infty(\vfi^{-1}) } (\vfi(v,y)+ \vfi(\tv, \ty))  \qquad \forall  t>0,\,\,v,y,\tv,\ty \in \R^d\,,
\ee
where $\omega, K$ depend  on $\ell_B, \ell_H, \gamma,$ and $\vfi$.
\end{theorem}

	\begin{proof} The strategy is to apply Theorem \ref{LTO} to $u(t+1)$  and reach the result by estimating $[u(1)]_{\theta}$ with Proposition \ref{hypocoer2}. {We first assume that $u_0$ is bounded, so that {$u(t)$} is not only bounded, but also  Lipschitz for $t>0$ (see Theorem \ref{teo-ex}). Next, we   apply Proposition \ref{hypocoer2}  (for $t=1$), which gives}, according to \rife{vfi3}, 
\begin{align*}
	& u(1,v,y)- u(1,\tv,\ty)   =  \\ 
	& =  \int_0^1 \left\{\nabla_v u(1,sv+(1-s)\tv, sy+(1-s)\ty)\cdot (v-\tv)+ \nabla_y u(1,sv+(1-s)\tv, sy+(1-s)\ty)\cdot(y-\ty)\right\}ds
	\\ & \qquad\qquad \leq   C\,    \|u_0\|_{L^\infty(\vfi^{-1})} \int_0^1 \vfi(sv+(1-s)\tv, sy+(1-s)\ty) ( |v-\tv|+ |y-\ty|)  \, ds 
	\\ &  \qquad\qquad  \leq \tilde C \,    \|u_0\|_{L^\infty(\vfi^{-1})}[\vfi(v,y)+ \vfi(\tv, \ty)]  ( |v-\tv|+ |y-\ty|)\,.
\end{align*}
We also notice that if $k_0$ satisfies $ (\partial_t + \cL)\vfi \geq - k_0$, then $\|u_0\|_{L^\infty(\vfi^{-1})} e^{k_0 t} \vfi$ is a supersolution of the equation since $\vfi\geq 1$. By comparison, we readily get 
\be\label{sabaudia}
	|u(t,v,y)| \leq  e^{k_0}\, \|u_0\|_{L^\infty(\vfi^{-1})} \vfi(v,y) \qquad \forall t\in (0,1]\,.
\ee
Thus, for a possibly different  constant $\tilde C$, we have
$$
	u(1,v,y)- u(1,\tv,\ty) \leq  \tilde C \,    \|u_0\|_{L^\infty(\vfi^{-1})}[\vfi(v,y)+ \vfi(\tv, \ty)]  \left( ( |v-\tv|+ |y-\ty|) \wedge 1\right)\,,
$$
which implies 
$$
	[u(1 )]_{1}\leq \tilde C \,    \|u_0\|_{L^\infty(\vfi^{-1})}\,.
$$
Finally, once we apply Theorem \ref{LTO} to $u(t+1)$, we get \rife{goal} for $t\geq 1$. However, up to increasing the value of $K$,  \rife{goal} also holds for $t\in (0,1]$ in light of \rife{sabaudia}.  { The case when $u_0$ is not bounded is obtained by approximation, as in Corollary \ref{visco-unbounded}.}
\end{proof}
\qed

Consequently, we deduce by duality our main result, which includes Theorem \ref{main} for the case that $\kappa=0$.

\begin{theorem}\label{main-new}  Assume that $H,B$ satisfy  \rife{Blip}, \rife{H1}-\rife{H2} and there exists a Lyapunov function $\vfi$ (cf. Definition \ref{def-Lyap}) satisfying conditions \rife{liplocvfi}, \rife{vfi2}, \rife{vfi3}. 
Then, there exist $\omega, K>0$ such that, for every initial data $m_{01}, m_{02}\in \cP_1(\R^{2d})$ for which  $\vfi\in L^1(dm_{0i})$, $i=1,2$, the corresponding solutions $m_1,m_2$ of \rife{FP} satisfy 
\be\label{output2}  
\| {m_1(t)}-{m_2(t)}\|_{TV_\vfi} \leq K \, e^{-\omega t} \, \|m_{01} - m_{02}\|_{TV_\vfi}, 
\ee
where $\|\cdot\|_{TV_\vfi}$ is the norm of the total variation weighted by $\vfi$.
\end{theorem}

\vskip1em
{\bf Proof of Theorem \ref{main-new} (and Theorem \ref{main}).}\quad 
We first observe that if $\vfi\in L^1(dm_{0i})$, $i=1,2$, then we also have $\vfi\in L^1(dm_{i}(t))$ for all $t>0$. This is established as in the proof of Theorem \ref{teo-ex}, where we have shown that if the initial measure has finite $k$ moments, then the corresponding solution enjoys the same property at later times. Setting $m_0\coloneqq m_{01}-m_{02}$ and $m\coloneqq m_1-m_2$, applying \rife{dual-form} for both solutions and subtracting, we obtain
$$
	\into \zeta\, dm(t) = \into {u(0)}\, d{ m}_0\,,
$$
for every $\zeta\in C_b(\R^{2d})$, where $u$ is the solution of the backward problem (B) in \rife{dual-form}.  By stability of viscosity solutions, standing on Corollary \ref{visco-unbounded}, and since $\vfi\in L^1(d|m|(t))$, the equality extends to the case $\zeta \in L^{\infty}(\vfi^{-1})$.  
Using that ${m}_0$ has zero-average, any constant can be added to {$u(0)$}. Then, it holds 
\be\label{Hai}
\begin{split}
	\into \zeta\, d m (t) &  \leq    \|m_0\|_{TV_\vfi}\, \inf_{c\in \R} \| {u(0)}+c\|_{L^\infty(\vfi^{-1})} 
	\\ & =  \|m_0\|_{TV_\vfi}\,\sup_{z, \tilde z} \frac{| u(0,z ) - u(0, \tilde z) |}{\vfi(z)+ \vfi(\tilde z)} \,
	\\ & \leq  \|m_0\|_{TV_\vfi}\,K\, e^{-\omega t}  \|\zeta\|_{L^\infty(\vfi^{-1}) }\,,
\end{split}
\ee
using Theorem \ref{LTO2}, where $z=(v,y)$, $\tilde z= (\tilde v, \ty)$, { and $\inf_{c\in \R} \| {u(0)}+c\|_{L^\infty(\vfi^{-1})} = \sup_{z, \tilde z} \frac{| u(0,z ) - u(0, \tilde z) |}{\vfi(z)+ \vfi(\tilde z)}$ follows from, e.g., \cite[Lemma 2.1]{HM}}. By arbitrariness of $\zeta$, this yields $\| \vfi\|_{L^1(d| {m|(t)} )}\leq  K\, e^{-\omega t} \|m_0\|_{TV_\vfi}$, which is the desired conclusion.
\qed

{ The result of Theorem \ref{main-new} can also be rephrased as the decay of zero-average solutions. Those solutions are meant in the sense of Corollary \ref{cor-ex} (possibly signed solutions). Hence, the proof of the result below works by duality,  exactly as in \rife{Hai}.

\begin{theorem}\label{zero-average} Assume the hypotheses of Theorem \ref{main-new} are valid. Let $m_0\in \cM(\R^{2d})$ be a measure with zero-average and such that $\vfi\in L^1(d|m_0|)$. Let $m$ be a solution of \rife{FP} in the sense of Corollary \ref{cor-ex}, i.e.,  $m\in C([0,T];\cM(\R^{2d})^*)$  and satisfies \rife{dual-form}. Then, we have
$$
\| {m(t)} \|_{TV_\vfi} \leq K \, e^{-\omega t} \, \|m_{0}  \|_{TV_\vfi}\,, \qquad \forall t>0\,.
$$
\end{theorem}
}

\subsection{Examples}\label{examples}

We conclude this section by giving a model class of examples, with a  specific Lyapunov function $\vfi$,  for  which the above theorem can be used. Precisely, we assume that
$$
	H(v,y)= v\,,\qquad B(v,y)= b(v,y)+ \nabla \Phi(y)\,,
$$
where $b$ is a Lipschitz function satisfying  
\be\label{bdiss}
	\exists \,\, \alpha,c_0>0\,: \quad b(v,y)\cdot v \geq \alpha \, |v|^2- c_0\,,\qquad b(v,y)\cdot y \geq -c_0 |y| (1+ |v|)\,,  \qquad \forall v,y\in \R^d\,, 
\ee
%
and $\Phi\in W^{2,\infty}(\R^d)$ fulfils 
%
\be\label{Phidiss}
	\exists \,\, \beta,c_1>0\,: \quad \nabla \Phi(y)\cdot y \geq \beta \, |y|^2- c_1 \qquad \forall  y\in \R^d \,.
\ee
Then, we can build a Lyapunov function  for \rife{pb}.

\begin{lemma}\label{lya} Suppose that \rife{bdiss} and \rife{Phidiss} hold. There exist $\vep, \de>0$ such that
$$
	\vfi:= \Phi(y) +  \frac12 \left( |v|^2 + 2\vep v\cdot y+ \de |y|^2\right)
$$
is a Lyapunov function for the operator $\cL$ satisfying \rife{liplocvfi}, \rife{vfi2}, and \rife{vfi3}.
\end{lemma}

\begin{proof} We compute and estimate, by virtue of \rife{bdiss} and \rife{Phidiss},
\begin{align*}
	\cL [\vfi] & = - d-{\kappa\de d}-x\cdot \nabla \Phi(y)  -  \vep |v|^2 - \de v\cdot y   +   b(v,y) \cdot v + \vep b(v,y) \cdot y  + \nabla \Phi(y)  \cdot (v+\vep y) 
	\\  & \geq ( \alpha-\vep)  |v|^2 - \de v\cdot y + \vep (\beta |y|^2- c_1) + \vep b(v,y) \cdot y -  { [d(1+\kappa\de)  +c_0]}
	\\ & \geq \bigg( \frac\alpha2-\vep\bigg)  |v|^2 + \vep \bigg(\frac\beta 2 |y|^2- c_1\bigg) - \frac{\de^2}{2\alpha} |y|^2  - \vep \frac{c_0^2}{2\beta} (1+ |v|)^2  -  { [d(1+\kappa\de)  +c_0]}
	\\ & \geq   \bigg( \frac\alpha2-\vep- \vep \frac{c_0^2}{\beta}\bigg)  |v|^2  + \vep  \bigg(\frac\beta 2- \frac{\de^2}{2\alpha\vep }\bigg)  |y|^2- \bigg(d{(1+\kappa\de)}  +c_0+ \vep c_1 +\vep \frac{c_0^2}{\beta}\bigg)\,.
\end{align*}
We now take $\de= \vep$ and choose $\vep<1$ sufficiently small (only depending on $\alpha, \beta,$ and $c_0$) so that 
$$
\frac\alpha2-\vep- \vep \frac{c_0^2}{\beta}\geq \frac\alpha 4\,,\qquad \frac\beta 2- \frac{\vep}{2\alpha }\geq \frac\beta4\,.
$$
Then, we get
\be\label{lbvfi} 
\cL[\vfi] \geq  \frac\alpha 4 |v|^2+ \vep \frac\beta4 |y|^2 - C_\vep\,,
\ee
and, being $\vep<1$, we also have $\de-\vep^2 = \vep-\vep^2>0$, which implies $\vfi\geq \Phi(y) + c_\vep ( |v|^2+ |y|^2)$.  Of course, we can assume $\Phi(y)\geq 1$ with no loss of generality (anyway, $\Phi$ is bounded below due to \rife{Phidiss}). In addition, we recall that $\nabla \Phi$ is Lipschitz, thus $\vfi \leq C (|v|^2+ |y|^2+1)$. Therefore, \rife{lbvfi} yields that $\vfi$ is a Lyapunov function, {see \rife{lyapu}}.  Moreover, {\rife{liplocvfi} is satisfied because $\nabla \Phi$ is Lipschitz continuous}. As for \rife{vfi2}, it is enough to observe that 
$$
\frac12 \cL (\vfi^2)= \vfi \cL (\vfi) - (|\nabla_v \vfi|^2+{\kappa|\nabla_y \vfi|^2}) 
$$
and the right-hand side diverges at infinity, due to \rife{lbvfi} and the quadratic growth of $\vfi$. Finally, \rife{vfi3} holds since $\vfi\in W^{2,\infty}(\R^{2d})$ (being the sum of a quadratic function and $\Phi\in W^{2,\infty}(\R^d)$) and hence semiconvex. Then, by just using convexity, \rife{vfi3} holds for $\vfi$ up to adding a quadratic function;  but this means that it holds for $\vfi$ itself because $\vfi$ dominates a quadratic function. 
\end{proof}
\qed

\begin{remark}
We observe that the above lemma admits  similar extensions to the case of possibly more  general functions $H(v)= \nabla \Psi(v)$, for uniformly convex $\Psi$, and to possibly superquadratic growth of the function $\Phi(y)$.  
\end{remark}

As a particular case,  we can now deduce from Theorem \ref{main} the exponential decay for the classical Kolmogorov equation \rife{KO}.

\begin{corollary} Assume that $b$ satisfies \rife{Blip} and \rife{bdiss} and $\Phi\in W^{2,\infty}(\R^d)$ fulfils \rife{Phidiss}.  Let $m_0\in \cM_2(\R^{2d})$  with $\int m_0=0$, and let $m$ be the solution (in the sense of Corollary \ref{cor-ex}) of \rife{KO} with $B(v,y)=b(v,y)+ \Phi(y)$. Then, there exist $\omega, K>0$ such that 
\be\label{FP-decay}  
{ \| {m(t)}\|_{\cM_2} \leq K \, e^{-\omega t} \, \|m_0\|_{\cM_2}}\,.
\ee
\end{corollary}

\begin{proof} Let  $\vfi$ be defined from Lemma \ref{lya}. At infinity, $\vfi$  is equivalent to $|v|^2+ |y|^2$, so \rife{FP-decay} can be  rephrased as
$$
\| {m(t)}\|_{TV_\vfi} \leq K \, e^{-\omega t} \, \|m_0\|_{TV_\vfi}\,,
$$
which is exactly the estimate of Theorem \ref{zero-average}.
\end{proof}

\section{Appendix}

In this appendix, we provide short proofs of the existence and uniqueness results for both problems \rife{FP} and \rife{pb}. We start with the latter one.

\begin{proposition}\label{viscos} Assume that \rife{Blip} and \rife{H1}-\rife{H2} hold, and $\kappa\geq 0$.  For any $u_0\in C_b(\R^{2d})$, there exists a unique viscosity solution $u\in   C_b([0,T] \times \R^{2d})$ of \rife{pb}.
\end{proposition}

\begin{proof} The existence is given through vanishing viscosity. Indeed, if $\kappa>0$ in \rife{pb}, then the existence of a unique (classical) solution is well known, due to the global Lipschitz conditions \rife{Blip} and \rife{H2} on the drift terms. Moreover, if $u_0$ is bounded and Lipschitz, then the solution is also bounded and Lipschitz, see, e.g., \cite[Section 3]{PoPr}. For the general case (possibly degenerate), we consider a Lipschitz and uniformly bounded sequence $u_{0n}$ of initial data, which converges to $u_0$ locally uniformly, and we take the solutions $u_n$ corresponding to $\kappa+\frac1n$. By the maximum principle, we have that {$u_n$} is uniformly bounded. In addition, being $u_n$ Lipschitz and smooth, according to Proposition \ref{hypocoer}, we achieve that {$\nabla_vu_n(t), \nabla_y u_n(t)$} are uniformly bounded for $t>0$. It follows that the sequence $u_n(t, \cdot)$ is equicontinuous. Now, applying \cite[Lemma 6.2]{Po} to $u(t+\cdot, \cdot)$, we deduce that, for every compact set $K\subset \R^{2d}$ and $\vep>0$, there exists $\de>0$ such that 
$\|{u_n(t+s )- u_n(t )} \|_{L^\infty(K)}\leq \vep$, for any $s\in (0,\de)$ and any $n$, where $\de$ only depends on the local bounds of $H,B,$ and $u_n$, and on the equicontinuity of $u_n(t )$. In particular, by Ascoli-Arzel\`a theorem, the sequence $u_n$ turns out to be relatively compact, in the uniform topology, in  $[t, T]\times K$, whatever $t>0$ and compact set $K$ are given. Selecting compact sets $K_j$ exhausting $\R^{2d}$ and a sequence $t_j \downarrow 0$, by a diagonal procedure, we can extract a  subsequence $u_{n_k}$ converging to  a continuous function $u$ locally uniformly in ${(0,T)}\times \R^{2d}$. Of course, from the previous estimates, $u$ is bounded and Lipschitz for $t>0$. Using again \cite[Lemma 6.2]{Po}, and in view of the local uniform convergence of $u_{0n}$ to $u_0$, we also deduce that $u$ is actually continuous up to $t=0$ and ${u(t)}\to u_0$ as $t\to 0$. Thus, we have $u\in   C_b([0,{T}] \times \R^{2d})$, and, by the stability of the viscosity formulation for local uniform convergence, we also get that $u$ is a viscosity solution of \rife{pb}.

To prove uniqueness, we use classical tools in the theory of viscosity solutions, so we only sketch the argument. First of all, we observe that, due to the linear growth of the drift terms $B$ and $H$, the function $e^{\lambda t}(|v|^2{+|y|^2}+ M)$ is a supersolution of the equation, up to choosing $\lambda$ and $M$ sufficiently large. Therefore, for any small $\vep$, $u_\vep:= u-\vep e^{\lambda t}(|v|^2{+|y|^2}+ M)$ and $v_\vep:= v+ \vep e^{\lambda t}(|v|^2{+|y|^2}+ M)$ are a strict subsolution, respectively supersolution, being $u_\vep \leq 0$ and $v_\vep \geq 0$ outside of a suitable compact set $K_\vep$. Thus, up to replacing $u,v$ with $u_\vep, v_\vep$, we are reduced to proving a comparison argument assuming that $\sup(u-v)$, if ever positive, would be attained in some compact set, where the classical doubling method for viscosity solutions allows to conclude. For the interested reader, we refer to \cite[Proposition 6.1]{Po} for a similar proof in a much more general context, with purely continuous data. 
\end{proof}
\qed

The next result shows that weak solutions of the Fokker-Planck equation  are in duality with the viscosity solutions constructed above.

\begin{lemma}\label{dualok}
		
Assume that \rife{Blip} and \rife{H1}-\rife{H2} hold, and $\kappa\geq 0$. Let $m_0\in \cM_1(\R^{2d})$, and let $m\in C([0,T];\cM(\R^{2d})^*)$ be a weak solution of \rife{FP} such that ${m(t)}\in \cM_1 (\R^{2d})$ and $ \sup\limits_{t\in (0,T)} \intd |z|\, d{|m|(t)}  <\infty$. {Then, $m$ satisfies \rife{dual-form}}.

\end{lemma}

\begin{proof}  In what follows, we denote $z:=(v,y)$ and we set $F=(B, -H)$ the drift term in $\R^{2d}$ for the equation satisfied by $m$.  We introduce   a sequence  $ \rho_\de(z)$ of standard  (compactly supported) symmetric mollifiers in $\R^{2d}$  and we set 
$$
m_\de(t,z)=   m(t)\star \rho_\de : = \int_{\mathbb R^{2d}}   m(t,w) \rho_\de(z-w)\, dw \,.
$$
We  also take  a sequence of 1-d mollifiers $\xi_\vep(t)$ such that supp$(\xi_\vep)\subset (-\vep,0)$, and we set 
$$
 m_{\de,\vep}:= \int_0^T \xi_\vep(s-t)\,m_\de(s)ds\,.
$$
One can easily verify (cf. \cite{P-BUMI}[Lemma 3.2]) that $\mdep$ satisfies the following equation
\be\label{eqconv}
\begin{cases}
\partial_t \mdep  -  \Delta_v \mdep - \kappa \Delta_y \mdep -{\rm div}_{(v,y)} ((m\, F)_{\de,\vep} )=- \xi_\vep(-t) [m_0\star\rho_\de]\,, & \hbox{in $(0,T)\times \R^{2d}$,}
\\
\mdep(0)= 0 & \hbox{in $ \R^{2d}$}
\end{cases}
\ee
where $(m\, F)_{\de,\vep}= (mF)_\de  \star\xi_\vep$. The fact that  $\mdep(0)=0$ is due to  the decentered choice of the time mollification; since $\xi_\vep$ is compactly supported in $(-\vep,0)$, it follows that $\mdep(t)$ vanishes for $t$ in a right neighborhood of zero and  we have
\be\label{tpicc}
\intd\mdep(t,z)\psi\, dz =  \intd\int_0^T (\xi_\vep(s-t)-\xi_\vep(s)) \,\psi\, m_\de(s,z)dsdz  \leq c_\vep \, t\, \|\psi\|_\infty\,,\quad \forall \psi\in C_b(\R^{2d})\,,\,\, \forall t>0
\ee
for some constant $c_\vep$ only depending on $\vep$. 

We now consider the vanishing approximation of the viscosity solution, namely the unique classical solution of the backward problem
$$
\begin{cases}
	-\partial_t u_n+ \cL[u_n]     =  \frac1n \Delta_yu_n &  \text{ in } ( 0,t ) \times \R^{2d}, \\ 
	u_n(t,v,y) = \zeta (v,y)\qquad  &  \text{ in }   \R^{2d}.
\end{cases} 
$$
In a first step, here one can consider that $\zeta$ is also a smooth Lipschitz function. In order to exploit the duality between $u_n$ and $\mdep$, we need to take an auxiliary cut-off function $\xi_R(z):= \xi (z/R)$ for some smooth function $\xi$ supported in $\{|z|\leq 2\}$ and such that $\xi\equiv 1$ in $\{|z|\leq 2\}$. Mutiplying the equation of $u_n$ by $\mdep\, \xi_R$ and integrating by parts (thanks to $\xi_R$), we get
\begin{align*}
 -\intd \zeta\, \xi_R\, d\mdep(t)+ \int_0^t\intd u_n\left[ \partial_t -  \Delta_v- (\kappa+ \frac1n) \Delta_y \right] (\mdep\, \xi_R) +  \int_0^t\intd F(z) \cdot\nabla_z u_n\, \mdep\, \xi_R = 0
\end{align*}
which yields, using the equation \rife{eqconv},
\begin{align*}
 & -\intd \zeta\, \xi_R\, d\mdep(t)   +  \int_0^t\intd F(z) \nabla_z u_n\, \mdep\, \xi_R - \int_0^t\intd (m\, F)_{\de,\vep} \nabla_z (u_n\xi_R)
 \\ 
 & \qquad =   \int_0^t\intd u_n\, \xi_\vep(-t) (m_0)_\de  \xi_R +  \frac1n \int_0^t\intd u_n\Delta_y\mdep \, \xi_R \\
 & \qquad + \int_0^t\intd u_n\mdep\, \left[ \Delta_v +(\kappa+ \frac1n) \Delta_y \right] (\xi_R) +  2u_n\, \left[\nabla_v \mdep\nabla_v \xi_R+ (\kappa+ \frac1n)\nabla_y \mdep\nabla_y \xi_R\right]
\end{align*}
As a first step, we get rid of the cut-off function by letting $R\to \infty$. This is justified by Lebesgue's theorem, using that $u_n$ is bounded and Lipschitz, and that 
$F(z) \mdep(z) \in L^1((0,T)\times \R^{2d})$\footnote{thanks to the Lipschitz character of $F$, and the fact that $\sup\limits_{t\in (0,T)} \intd |z|\, d{|m|(t)}  <\infty$;  indeed, we have $|F(z) \mdep(z)|\leq C [((1+ |z|)|m|)\star\rho_\de]\star\xi_\vep$}. In particular, all the terms with derivatives of $\xi_R$ vanish as $R\to \infty$. Then we get 
\be\label{dualn}
\begin{split}
 & -\intd \zeta\,  d\mdep(t)   +  \int_0^t\intd F(z) \nabla_z u_n\, \mdep\,  - \int_0^t\intd (m\, F)_{\de,\vep} \nabla_z  u_n 
 \\ 
 & \qquad =   \int_0^t\intd u_n\, \xi_\vep(-t) (m_0)_\de  +  \frac1n \int_0^t\intd u_n\Delta_y\mdep \, 
 \end{split}
 \ee
Now we estimate, using the Lipschitz regularity of $F$ and Proposition \ref{hypocoer} for $\nabla u_n$,
\be\label{commu}
\begin{split}
 &
\qquad  \left |\int_0^t\intd F(z) \nabla_z u_n\, \mdep\,  - \int_0^t\intd (m\, F)_{\de,\vep} \nabla_z  u_n \right|  \\ 
 & = \int_0^t \int_0^t \xi_\vep(s-\tau)\intd \intd \nabla_zu_n(\tau, z) [ F(z)-F(w)]   \rho_\de(z-w)\,  dm(s,w)\, dz d\tau
 \\ & \leq c\, \de\,  \int_0^t  \intd  \tau^{-\frac32} m_{\de,\vep}(\tau,z)\, dz d\tau \leq \tilde c_\vep\, \, \de
 \end{split}
 \ee 
where last inequality comes from \rife{tpicc}.  Now we pass to the limit in \rife{dualn}. Here  we first let $n\to \infty$ (so last term vanishes since $\mdep$ is smooth), using the uniform convergence of $u_n$. Secondly we let $\de \to 0$, at fixed $\vep$, using estimate \rife{commu}. We obtain so far
$$
-\intd \zeta\,  dm_\vep(t)    =   \int_0^t\intd u(t) \, \xi_\vep(-t)  dm_0 
 $$
 and finally we let $\vep \to 0$, using the continuity of $u(t)$ (recall that $u(t)\to u(0)$ locally uniformly), the continuity of  $m(t)$ and the fact that $m$ has finite first moment. We conclude that \rife{dual-form} holds true, for any $\zeta\in C_b(\R^{2d})$ and for the (unique) viscosity solution corresponding to $\zeta$.
\end{proof}
\qed

In light of the above lemma, we conclude with the proof of Theorem \ref{teo-ex}.
 \vskip0.3em
 
{\bf Proof of Theorem \ref{teo-ex}.}\quad Item (i) was  proved in Proposition \ref{viscos}. As for the Fokker-Planck equation, the uniqueness is a straightforward consequence of Lemma \ref{dualok} whenever $ m(t)\in \cM_1 (\R^{2d})$ for $t\in [0,T]$  and $ \sup\limits_{t\in (0,T)} \intd |z|\, d{|m|(t)}  <\infty$. Indeed, if $m_0=0$, then {\rife{dual-form}} implies $\intd \zeta d{m(t)}=0$ for any $t$. Here, any $\zeta \in C_b(\R^{2d})$ is allowed, due to the solvability of the adjoint problem, which is guaranteed by item (i). Hence, ${m(t)}=0$.

If $m\in C([0,T];\cP(\R^{2d}))$, then we do not need to assume a priori that $m$ has finite first moment, and we now show this fact. Simply using the weak formulation and the 
time continuity of $m$, we have
\be\label{abov}
\intd \phi(t, v,y)\, d{m(t)}+ 
 \int_0^t\intd  \left( -\partial_t \phi+ \cL[\phi] \right) \,d{m(\tau)} \,d\tau = \intd {\phi(0,v,y)}\, dm_0  
 \qquad\forall \,  \phi\in C^{1,2}_c(\overline Q_t)\,,\,\, t>0\,.
\ee
We wish to remove the condition that $\phi$ has compact support in $\R^{2d}$.  We set  $\langle z\rangle:= \sqrt{1+|z|^2}$ and take a real-valued, cut-off function $\xi(s)$, {supported in $\{|s|\leq 2\}$ and such that $\xi\equiv 1$ in $\{|s|\leq 1\}$}. Then, we replace $\phi$ with  $\phi\,  \xi_R$ as a test function in \rife{abov}, where $\xi_R(z):= \xi\left(\frac{\langle z\rangle}R\right)$ and now $\phi\in C^{1,2}(\overline Q_t) $ is assumed to be nonnegative but does not need to have compact support in $\R^{2d}$.  We achieve 
$$
\intd  \phi \xi_R \, d{m(t)}+ 
 \int_0^t\intd  \left( -\partial_t \phi+ \cL[\phi] \right) \xi_R\,d{m(\tau)} \,d\tau   =  \intd {\phi(0,v,y)}\, \xi_R\, dm_0+ o_R(1)\,,
 $$
where $o_R(1)$ is a quantity converging to zero as $R\to \infty$, since it involves the derivatives of $\xi_R$, which vanish uniformly as $R\to \infty$, with the fact that $\phi\in C^{1,2}(\overline Q_t)$, and $m$ has finite mass. Letting $R\to \infty$ leads to 
\be\label{gene}
 \intd  \phi   \, d{m(t)}+ 
 \int_0^t\intd  \left( -\partial_t \phi+ \cL[\phi] \right) \,d{m(\tau)} \,d\tau   =  \intd {\phi(0,v,y)}\,  \, dm_0\,.
 \ee
Notice that  if $\phi= e^{\lambda(T-t)}( \langle z\rangle^k   + M)$, then one can verify that $(-\partial_t + \cL)[\phi]\geq c\, \langle z\rangle^k \geq 0$, by choosing $\lambda$ and $M$ large enough, and the same holds replacing $\langle z\rangle^k$ with $\langle z\rangle^k \wedge n$, for any $n>0$. Therefore, by virtue of \rife{gene} with $\phi= e^{\lambda({T-t})}( (\langle z\rangle^k \wedge n)   + M)$, it holds
$$
 \intd  (\langle z\rangle^k \wedge n)  \, d{m(t)}\leq C\, e^{\lambda t}   \intd \langle z\rangle^k \,  \, dm_0\,.
$$
This yields $ \intd |z|^k d{m(t)}  \leq C_T\, \intd |z|^k\, dm_0$ for all $t\in (0,T)$. In particular, we have established that
 if $m_0$ has $k$ finite moments, then any {\it nonnegative} solution ${m(t)}$ has $k$ finite and bounded moments in $(0,T)$. This concludes {the proof of} the uniqueness statement of Theorem \ref{teo-ex}. We also note that, with the same estimate as above, one can  prove that $m$ is continuous from $[0,T]$ into $\cP_1(\R^{2d})$.

We are only left to show the existence of at least one solution,  which we  construct  by vanishing viscosity limit. Assume that $m_0\in \cP{(\R^{2d})}$ has $k$ finite moments. Let $m_n$ be  the unique (classical) solution of \rife{FP} when replacing $\kappa$ with $\kappa+\frac1n$.  By the above estimate, we know that 
\be\label{mome}
 \intd \langle z\rangle^k\,{dm_n(t)}  \leq C_T\, \intd \langle z\rangle^k\, dm_0\leq C\,,
\ee
which gives uniform bounds for the {$k$th} moments of $m_n$. With standard arguments, see, e.g., \cite[Lemma 6.3]{Po}, one can prove that $m_n$ is relatively compact in $C([0,T];\cP(\R^{2d}))$, and, up to a subsequence, it converges to a weak solution $m$. In addition, ${m(t)}\in \cM_k$ for every $t\geq 0$.  For a general signed measure $m_0$, of course, we reason by splitting $m_0= m_0^+-m_0^-$ and obtaining a solution for both the negative and positive parts.
\qed

\vskip1em

{\bf Proof of Corollary \ref{cor-ex}.}\quad  If $m_0\in \cP({\R^{2d}})\cap\cM_k{(\R^{2d})}$, then it can be approximated by  a sequence $m_{0n}$ of initial data belonging to $\cP_1{(\R^{2d})}$ that have $k$ uniformly bounded moments.  Reasoning as in the previous proof, the corresponding sequence of solutions, $m_n$, satisfies a uniform  estimate such as \rife{mome}, and it is relatively compact in $C([0,T];\cP(\R^{2d}))$. A subsequence thus converges to some $m$, and since equality {\rife{dual-form}} holds for $m_n$, passing to the limit we get it for $m$ as well. This result extends to any nonnegative measure $m$, and, for a general $m_0$, one can reason by splitting  $m_0= m_0^+-m_0^-$ in order to construct a solution. Finally, the uniqueness  is a straightforward consequence of {\rife{dual-form}} because the adjoint problem is solvable for all $\zeta\in C_b(\R^{2d})$. The preservation of positivity follows from the adjoint problem, and the {conservation} of mass is obvious being $u=1$ the only viscosity solution corresponding to $\zeta=1$.
\qed

\end{document}